\documentclass[11pt,reqno]{amsart}

\usepackage[x11names]{xcolor}
\usepackage[linkcolor=Blue2,citecolor=black,colorlinks]{hyperref}
\usepackage{color, bm, amscd, tikz-cd}

\usepackage{mathdots}
\usepackage{amsmath,amsthm,amssymb, mathrsfs, bbm}

\allowdisplaybreaks

\theoremstyle{plain}
\newtheorem{theorem}{Theorem}[section]
\newtheorem{corollary}[theorem]{Corollary}

\newtheorem{proposition}[theorem]{Proposition}

\theoremstyle{definition}
\newtheorem{definition}[theorem]{Definition}
\newtheorem{example}[theorem]{Example}

\newtheorem{remark}[theorem]{Remark}

\numberwithin{equation}{section}

\newcommand{\R}{\mathbb{R}}
\newcommand{\C}{\mathbb{C}}

\newcommand{\D}{\mathbb{D}}
\newcommand{\T}{\mathbb{T}}

\newcommand{\bydef}{\stackrel{\rm def}{=}}
\usepackage{mismath}

\begin{document}
\title[Approximation]{The multivariate Herglotz-Nevanlinna class:  rational approximation}

\author[M. Bhowmik]{Mainak Bhowmik }
\address{Indian Institute of Science, Bangalore, India }
\email{\tt mainak.bhowmik943@gmail.com, mainakb@iisc.ac.in}

\author[M. Putinar]{Mihai Putinar}
\address{University of California at Santa Barbara, CA} 

\email{\tt mputinar@math.ucsb.edu}



 

\begin{abstract} 

 We return to Takagi's variational principle, generalized after forty years to two complex variables by Pfister. Both isolating some extremal rational functions associated to a bounded holomorphic function
in the unit disk, respectively the bidisk. The rational inner functions arising from the Takagi-Pfister skew eigenvectors lead to a Pad\'e type approximation scheme. For these rational functions, we prove a Montessus de Ballore type convergence theorem, on the polydisk in any complex dimension.
On the natural and more restrictive class of Agler holomorphic functions with non-negative real part, we show that Cayley rational inner functions match any finite section of the Taylor expansion at a prescribed point. We derive from the Hilbert space proof that the finite section coefficient set of Taylor series of the Agler functions in the Herglotz-Nevanlinna setting is semi-algbraic. The pole distribution of the Takagi-Pfister interpolation sequence is identified as a main open question on the subject.

\end{abstract}

\maketitle

\section{Introduction} 

This is a second article devoted to the multivariate Herglotz-Nevanlinna class of holomorphic functions. Below we focus on rational approximation, constructive as much as possible, of bounded holomorphic functions depending on several complex variables.
The polydisk as a domain of definition for these functions is privileged, due to the simple Schwarz reflection operation carried by the torus, its distinguished boundary.

From the rich, rather technical and still in progress realm of multivariate bounded holomorphic interpolation we touch two themes: the abundance of rational inner functions
in the polydisk and a variational principle attached to the Hankel form induced by a bounded holomorphic symbol. The later is exploited as a gate to a multivariate Pad\'e approximation scheme.
We closely follow two classical works: Takagi's century old article \cite{Takagi} and Pfister's doctoral dissertation of 1960 \cite{Pfister-thesis, Pfister}. Takagi's work timely complemented landmark contributions of
Schur, Carath\'eodory and Fej\'er  to the 1D bounded holomorphic interpolation theory, while Pfister's thesis is the first notable attempt to explore the multivariable counterpart
of these advances.
We return to Takagi and Pfister from our days perspective of accumulated innovation, such as the Hilbert space realization of structured holomorphic functions.

Let $f$ be  a holomorphic function defined in the polydisk $\D^d, d \geq 1,$ subject to a uniform bound condition, or a non-negative real part assumption. The Taylor series expansion
at the origin is
$$ f(\boldsymbol{z}) = \sum_{\boldsymbol{\alpha} \in \mathbb{N}_0^d} c_{\boldsymbol{\alpha}}(f) \boldsymbol{z}^{\boldsymbol{\alpha}}.$$
We explore below the interplay between a finite section of the coefficients $(c_{\boldsymbol{\alpha}}(f))$ and the imposed restrictions on $f$, with a final aim at approximating $f$ by rational functions matching the Taylor coefficient data.

We owe to Jim Agler \cite{Agler-Adv.Op.Th} for the discovery of the natural framework for doing multivariate bounded holomorphic interpolation, as closely as possible to the classical single variable case.
In this context we prove two results: a finite section of the Taylor series of a Herglotz-Nevanlinna function belonging to the Agler's class can be matched exactly by a Cayley rational inner function
of degree depending linearly on the size of the prescribed finite section (see Theorem \ref{Thm:Inner-solution}). And second, we infer that such a finite section of coefficients of all Taylor series (at the origin) of functions in Agler's class is semi-algebraic (see Theorem \ref{semialgebraic}).
In other terms, finite determinateness of different sorts is embedded into the coefficients of Agler functions. Definitions and full details are given in the sequel.

In the second part of the article, we propose the following Pad\'e type approximation scheme for any holomorphic function $f: \D^d \longrightarrow \overline{\D}$. We seek, in the spirit of Takagi and Pfister,  the optimal value
$\sigma_{2\boldsymbol{n}}$ 
of the Hankel form supported on the unit torus $\T^d$:
$$ \max \Re \int_{\T^d} f(\boldsymbol{\zeta}) q(\boldsymbol{\zeta})^2 dm(\boldsymbol{\zeta}), \ \ \| q \|_2 = 1,$$
where $q(\boldsymbol{\zeta})$ is a trigonometric polynomial of multi-degree $\boldsymbol{n}$. Assume this value is attained at an element $q_{2\boldsymbol{n}}$.
Then, under a weak-topology assumption on the family $(q_{2 \boldsymbol{n}})$, there exists an analytic hypersurface $X$ of $\D^d$, with the property that
a subsequence of the rational functions $ \sigma_{2\boldsymbol{n}} \frac{q_{2 \boldsymbol{n}}}{\overline{q_2 \boldsymbol{n}}}$ converges uniformly to $f$
on compact subsets of $\D^d \setminus X$ (see Theorem \ref{weak-conv}). This result aligns with the Montessus de Ballore type theorems studied in the multivariate setting \cite{Cuyt, Guillaume-Huard}.
The occurrence of a Hankel form in the approximation process is not unexpected, having the single variable analysis due to Adamjan-Arov-Krein as a comparison basis,
see for instance \cite{Peller-book}.

The contents are the following. Section \ref{Prelim} is devoted to preliminaries. In Section \ref{Existence of Cayley rational inner solution} we prove that Cayley rational inner functions in the polydisk $\D^d$ exactly match a finite section of a normalized Agler function in the Herglotz-Nevanlinna setting. The degree bound of this exact rational interpolant is linear with respect to size of the coefficient data. In Section \ref{The coefficient problem} we prove that the coefficient set of a prescribed finite section of normalized Agler functions in the polydisk is semialgebraic. A couple of low degree examples are included, reflecting the invariance of the coefficient convex body under the automorphism group of the right half plane. Section \ref{Rational approximation} introduces the Takagi-Pfister approximation scheme and relates it to a classical theme of Hankel form approximation. The inherent challenges of the multivariate analogue of a well understood single complex variable setting are analyzed. A couple of open questions and comments around them are included in Section \ref{Final remarks and open questions}.
\bigskip

{\bf Acknowledgements.} The first author was supported by the Prime Minister’s Research Fellowship PMRF-21-1274.03.The second author was partially supported by a Simons Collaboration Grant for Mathematicians.

\section{Herglotz-Nevanlinna class in the polydisk} \label{Prelim}

This section aims to recall some terminology and collect a series of known ingredients that will appear in the body of this note.
The excellent monograph \cite{Agler-McCarthy-Young} offers complete proofs and pertinent comments on these facts, all structured
in a natural framework.

\begin{definition} \label{Def: HN}
A holomorphic function $f: \D^d \rightarrow \mathbb{C}$ is said to be a 
\begin{enumerate}
\item[(i)] {\em Schur function}, provided $|f(\boldsymbol{z})| \leq 1$ for every $\boldsymbol{z} \in \D^d$;

\item[(ii)] {\em Herglotz-Nevanlinna function}, provided $\Re f(\boldsymbol{z}) \geq 0$ for every $\boldsymbol{z} ~\in~ \D^d$.
\end{enumerate}
\end{definition}
 The set of all Schur functions and the set of all Herglotz-Nevanlinna functions on $\D^d$ are denoted by $\mathcal{S}(\D^d)$ and $\mathcal{H}(\D^d)$, respectively. 
These two classes are in bijective correspondence via a Cayley transform $ w = \frac{z-1}{z+1}$ mapping the right half plane to the unit disk $\D$.

\subsection{Representation of Herglotz-Nevanlinna functions}
The celebrated Riesz-Herglotz representation theorem associates an element $\varphi \in \mathcal{H}(\D)$ to a unique positive regular Borel measure $\mu$ on the unit circle $\mathbb{T}$ via the integral representation:
\begin{align*}
\varphi(z)= i \Im \varphi(0) + \int_{\mathbb{T}} \frac{\xi + z}{\xi -z} d\mu(\xi) \, \, \text{ for } z\in \D.
\end{align*}
A counterpart of Riesz-Herglotz represenetation in the polydisk $\D^d$ was obtained by Kor\'anyi and Puk\'anszky \cite{KP}. It can be described as follows: for $\varphi$ in $\mathcal{H}(\D^d)$, there exists a unique positive regular Borel measure $\nu$ on the distinguished boundary $\mathbb{T}^d$ of $\D^d$ such that
\begin{align} \label{Eq:KP-rep}
			\varphi(\boldsymbol{z})= i\Im \varphi(0)+ \int_{\mathbb{T}^d}  \left(  \frac{2}{\prod_{j=1}^{d} (1 - z_j \bar{\xi}_j)}-1  \right) d\nu(\boldsymbol{\xi}),
\end{align}
where $\nu$ satisfies the moment conditions,
	\begin{equation} \label{Eq:Pluri-measure} 
\widehat{\nu}(n_1, \dots, n_d)= \int_{\mathbb T^d} \bar{\xi}_1^{n_1} \cdots \bar{\xi}_d^{n_d} d\nu(\boldsymbol{\xi}) = 0,
\end{equation} 
unless $n_j \geq 0$ for all $j=1, \cdots, d$ or  $n_j \leq 0$ for all $j=1, \cdots, d$.

As a matter of fact, for a Herglotz-Nevanlinna function $\varphi$, the radial limit 
$$
\lim_{r \to 1^-} \Re \varphi (r \boldsymbol{\xi})
$$
exists $\Theta_d$-a.e. on $\mathbb{T}^d$ and it is the density of the absolute continuous part of $\nu$ with respect to $\Theta_d$ (also sometimes denoted by $\Theta$) which is the normalized Haar measure on $\mathbb{T}^d$.
The integral representation above connects the Herglotz-Nevanlinna functions with the Carathéodory-Fejér interpolation and the moment problem naturally. From this perspective a superresolution phenomena has been studied by the authors in \cite{BP}.

Most of the challenges of multivariate bounded holomorphic interpolation (in the polydisk or in more general classical domains) reside in the rigidity of the moment vanishing conditions (\ref{Eq:Pluri-measure}). 
For instance not allowing a canonical approximation of these measures by simple ones, such as finitely many point masses.

\subsection{Agler class} \label{Sub: Agler class}
The main innovation came from Jim Agler \cite{Agler-Adv.Op.Th}, by restricting the set of bounded holomorphic functions to those well fitted for a commuting
tuple of linear bounded operators. Specifically, let $\varphi$ be a holomorphic function defined on $\D^d$ with Taylor series at the origin:
$$
\varphi(\boldsymbol{z})= \sum_{\boldsymbol{\alpha} \in \mathbb{N}_0^d} c_{\boldsymbol{\alpha}}(\varphi) \boldsymbol{z}^{\boldsymbol{\alpha}}.
$$ 
Let $\mathcal{E}$ be a fixed infinite dimensional, complex separable Hilbert space. For a $d$-tuple of commuting strict contractions $\underline{T}=(T_1, \dots, T_d)$ on $\mathcal{E}$ (i.e., $\|T_j\| < 1$ for all $j=1, \dots, d$), we define 
$$
\varphi(\underline{T}) = \sum_{\boldsymbol{\alpha} \in \mathbb{N}_0^d} c_{\boldsymbol{\alpha}}(\varphi)  \underline{T}^{\boldsymbol{\alpha}},
$$
where $\boldsymbol{\alpha} = (\alpha_1, \dots, \alpha_d)$, $\underline{T}^{\boldsymbol{\alpha}}= T_1^{\alpha_1} \dots T_d^{\alpha_d}$ and the above series converges in the operator norm topology. One can define a norm for the function $\varphi$ as follows:
\begin{align} \label{Agler-norm}
	\| \varphi \|_{\operatorname{Agler}} = \sup \{ \| \varphi(\underline{T})\| :  \underline{T} \text{ is a } d\text{-tuple of commuting strict contractions}\}.
\end{align}
Let $\mathcal{F}_d$ be the algebra of all holomorphic functions on $\D^d$ for which the Agler-norm \eqref{Agler-norm} is finite. The closed Agler-norm unit ball of $\mathcal{F}_d$ is called as the the {\em Agler} class, $\mathcal{A}\mathcal{S}(\D^d)$ in the setting of Schur function.

We say that the function $\varphi$ belongs to the {\em Agler} class, $\mathcal{A}\mathcal{H}(\D^d)$ in the Herglotz-Nevanlinna setting if 
\begin{align}\label{Eq:A-H-class}
\Re \varphi(\underline{T}) = \frac{1}{2} (\varphi(\underline{T}) + \varphi(\underline{T})^*) \geq 0,
\end{align}
for every $d$-tuple of commuting strict contractions $\underline{T}$ on $\mathcal{E}$. 



Note that, if $T_j = z_j I_{\mathcal{E}}$ for $1\leq j\leq d$, then for $\varphi $ in $\mathcal{A}\mathcal{S}(\D^d)$ (respectively, in $\mathcal{A}\mathcal{H}(\D^d$)) we have $|\varphi(\boldsymbol{z})| \leq 1$ (respectively, $ \Re \varphi(\boldsymbol{z}) \geq 0$).
 The converse holds for $d=1$ by von Neumann's inequality and for $d=2$ by a deep theorem due to And\"o. However, for $d\geq 3$, $\mathcal{A}\mathcal{H}(\D^d) \subsetneq \mathcal{H}(\D^d)$ and $\mathcal{A}\mathcal{S}(\D^d) \subsetneq \mathcal{S}(\D^d)$.

The functions $\varphi$ belonging to $\mathcal{A}\mathcal{H}(\D^d)$ are characterized by a special, weighted
positive definite kernel decomposition of the function $\varphi(\boldsymbol{z}) + \overline{\varphi(\boldsymbol{w})}$, or equivalently, by a ``resolvent realization" of a linear pencil of operators. The later touching linear control theory. We refer to \cite{Agler-McCarthy-Young} for details.

\subsection{Carath\'eodory-Fej\'er interpolation}

Contractive holomorphic functions in the polydisk are determined by the full Taylor expansion at a fixed point (usually taken to be the origin). A series of natural questions arise: can one characterize the set of coefficients of these functions? Given a finite section of the coefficients, is it possible to construct a rational function which matches the given partial Taylor series? What is the rate of convergence of this sequence of rational (formal) approximations?

The answers of all these queries in the single variable case are classical, due to the works of Schur, Carath\'eodory and Fej\'er to name only the discoverers; all works remaining at the level of algebra or function theory. But there is much more: some inherent matrix positivity was embedded into the original proofs. These structured matrices (called by the names of Toeplitz or Hankel) have shaped over an entire century a rich Hilbert space operator theory construct whose ramifications continue to enchant and puzzle \cite{FF}.

Throughout this note we consider the entrywise partial order on $\mathbb{N}_0^d $:
$$
\boldsymbol{\beta} \leq \boldsymbol{n} \ \text{ if and only if } \beta_j \leq n_j \  \text{ for every } 1\leq j \leq d.
$$
For a fixed multi-index $\boldsymbol{n}$, $\Lambda_{\boldsymbol{n}} =\{ \boldsymbol{\beta} \in \mathbb{N}_0^d: \boldsymbol{\beta} \leq \boldsymbol{n}\}$.
The circulating terminology is the following. Given a set of scalars $\{a_{\boldsymbol{\alpha}}: \boldsymbol{\alpha} \in \Lambda_{\boldsymbol{n}} \}$ (a.k.a. the interpolation data), the {\em Carath\'eodory-Fej\'er interpolation problem} on $\D^d$ asks for conditions assuring  the existence of a function $\varphi$ in $\mathcal{A}\mathcal{H}(\D^d)$ or $\mathcal{H}(\D^d)$ with Taylor series
\begin{align*}
f(\boldsymbol{z}) = \sum_{\boldsymbol{\alpha} \in \mathbb{N}_0^d} c_{\boldsymbol{\alpha}}(\varphi) \boldsymbol{z}^{\boldsymbol{\alpha}} \, \text{ for } \boldsymbol{z} \in \D^d 
\end{align*}
matching the prescribed set of coefficients:  $c_{\boldsymbol{\alpha}}(f) = a_{\boldsymbol{\alpha}}$ for $\boldsymbol{\alpha} \in \Lambda_{\boldsymbol{n}}$.

In one dimension, contractivity of a Toeplitz matrix formed by the given data determines the solvability of the interpolation problem in the Schur class, while the positivity of a Hankel matrix with entries from the given data ensures solvability in the Herglotz-Nevanlinna class. See \cite{FF} for a detailed treatment. 

In higher dimension ($d>1$), all known solvability criteria of this interpolation problem have a non-constructive component \cite{Dautov-Khud, Eschmeier-Patton-Putinar, Woerdeman}. 
Restricting the range of functions to Agler's class is instrumental. In this context, we refer to Woerdeman's article \cite[Theorem 2.1]{Woerdeman} which establishes necessary and sufficient conditions for the
existence of a solution to the Carath\'eodory-Fej\'er interpolation problem in the Agler setting. We will exploit this advance in a subsequent section.

\subsection{Change of Taylor coefficients under Cayley transforms} \label{Subsection: Change of Taylor coeff.}
For any $f\in \mathcal{A}\mathcal{S}(\D^d)$ which is not identically $1$, define $g = (1+f)(1-f)^{-1}$. Then $g \in \mathcal{A}\mathcal{H}(\D^d)$. 
We note that for each $\boldsymbol{\alpha} \in \mathbb{N}_0^d$, the Taylor coeffcient $c_{\boldsymbol{\alpha}}(g)$ depends on the Taylor coefficients of $f$ in a universal manner i.e., there is a rational function $M_{\boldsymbol{\alpha}}$ depending on $|\Lambda_{\boldsymbol{\alpha}}|$ complex variables (which are $c_{\boldsymbol{\beta}}(f)$ for $0 \leq \boldsymbol{\beta} \leq \boldsymbol{\alpha}$) such that 
$$c_{\boldsymbol{\alpha}} (g) = M_{\boldsymbol{\alpha}}\left( (c_{\boldsymbol{\beta}}(f))_{\boldsymbol{\beta} \leq \boldsymbol{\alpha}} \right).$$ In fact, the denominator of the rational function $M_{\boldsymbol{\alpha}}$ is some non-negative integer power of $(1-c_{\boldsymbol{0}}(f))$. 
For example, when $d=2$ and $\boldsymbol{\alpha}=(1, 1)$,
$$
c_{\boldsymbol{\alpha}} (g) = \frac{\partial^2 g}{\partial z_1 \partial z_2} (0) =2 \frac{(1-c_{(0,0)}(f)) c_{(1,1)}(f) + 2 c_{(1,0)}(f) c_{(0,1)}(f)}{(1-c_{(0,0)}(f))^3}.
$$

In a similar fashion, for any $\varphi \in \mathcal{A}\mathcal{H}(\D^d)$, if we take $g=(\varphi -1)(\varphi +1)^{-1}$, then there are universal rational functions $Q_{\boldsymbol{\alpha}}$ depending on $|\Lambda_{\boldsymbol{\alpha}}|$ complex variables (which are $c_{\boldsymbol{\beta}}(\varphi)$ for $0 \leq \boldsymbol{\beta} \leq \boldsymbol{\alpha}$) such that 
$$c_{\boldsymbol{\alpha}} (g) = Q_{\boldsymbol{\alpha}}\left( (c_{\boldsymbol{\beta}}(\varphi))_{\boldsymbol{\beta} \leq \boldsymbol{\alpha}} \right).$$
In this case, the denominator of $Q_{\boldsymbol{\alpha}}$ is a non-negative power integer power of $(1+ c_{\boldsymbol{0}}(\varphi))$. 

\subsection{Cayley rational inner functions}

 Rational inner functions are desirable in interpolation problem since they are finitely determined (by a truncated power series expansion or finitely many point evaluations). 
 We have exploited this feature in our preceding article \cite{BP}. Now we focus on approximation by particular rational functions. Among these, the inner rational functions
 on the polydisk stand aside.

\begin{definition}
A bounded holomorphic function $f$ on $\D^d$ is said to be {\em rational inner} provided $f$ is rational with poles off $\D^d$ and is unimodular $\Theta_d-$a.e. on $\mathbb{T}^d$.
\end{definition}

For instance finite Blaschke products 
$$ \prod_{j=1}^m \frac{z-\lambda_j}{1- \overline{\lambda_j}z}, \ |\lambda_j| <1, \ 1 \leq j \leq m,$$
are rational inner functions on $\D$ and it can be shown that, up to a unimodular multiplicative factor, these are all rational inner functions in $\D$. 

Any rational inner function on $\D^d$ is of the form: 
$$  \boldsymbol{z}^{\boldsymbol{m}} \frac{p^*(\boldsymbol{z})}{p(\boldsymbol{z})} $$
for some polynomial $p$ in $\mathbb{C}[\boldsymbol{z}]$ of multi-degree $\boldsymbol{n}=(n_1, \dots, n_d)$ (say) such that $p$ is non-vanishing in $\D^d$. Here the polynomial $p^*$, called as the {\em reflection} of $p$, is defined as follows
$$
p^*(\boldsymbol{z}) = \boldsymbol{z}^{\boldsymbol{n}} \overline{p \left(\frac{1}{\overline{\boldsymbol{z}}} \right)}.
$$
The multi-degree $\boldsymbol{n}=(n_1, \dots, n_d)$ of a polynomial $p$ is defined by $n_j = \deg_{z_j} p, \ \ 1 \leq j \leq d.$ Note that, the coefficient of $\boldsymbol{z}^{\boldsymbol{n}}$ can be null. See \cite[Satz 2]{Pfister} for more details.
For an authoritative account of rational inner functions in the polydisk we refer to the recent survey article by Greg Knese \cite{Knese}. 

In the Herglotz-Nevanlinna setting the {\em Cayley rational inner functions} are analogous to the rational inner functions in the Schur class.
\begin{definition}
A holomorphic function $\varphi$ with non-negative real part on $\D^d$ is said to be Cayley rational inner if $\varphi$ is rational with poles off $\D^d$ and $\Re \varphi$ is zero $\Theta_d-$ a.e. on $\mathbb{T}^d$.
\end{definition}
It is easy to check that a Cayley rational inner function $\varphi$ is the Cayley transform $(z \mapsto (1+z)/(1-z))$ of the rational inner function, $f = (\varphi -1)/(\varphi +1)$. 

\section{Existence of Cayley rational inner solution} \label{Existence of Cayley rational inner solution}

We exploit below the full force of Agler's theory to show that there are sufficiently many inner rational functions in the polydisk to solve the Carath\'eodory-Fej\'er truncated interpolation problem with respect to Agler's class.
There is a price to pay for this existence result: the degree of the rational interpolant is proportional, but bigger, than the cardinality of the prescribed data.
\begin{theorem} \label{Thm:Inner-solution}
	Let $\boldsymbol{n} $ be in $\mathbb{N}_0^d$. Every interpolation data $\{c_{\boldsymbol{\beta}}:  \boldsymbol{\beta} \in \Lambda_{\boldsymbol{n}} , c_{\boldsymbol{0}}=1 \}$  which is solvable by a function in $\mathcal{A}\mathcal{H}(\D^d)$ has a Cayley rational inner solution of multi-degree at most $( | \Lambda_{\boldsymbol{n}}| d, \dots,  | \Lambda_{\boldsymbol{n}}| d)$.
\end{theorem}

For notational simplicity, we denote $ \Lambda_{\boldsymbol{n}} $ by $\Lambda$. Let $\mathbb{C}^{|\Lambda|}$ be the Hilbert space of $|\Lambda|$-tuples $\boldsymbol{\xi} = (\xi_{\boldsymbol{\beta}})_{\boldsymbol{\beta} \in \Lambda}$ with $\| \boldsymbol{\xi} \|^2 = \sum_{\boldsymbol{\beta} \in \Lambda} |\xi_{\boldsymbol{\beta}}|^2$. We prefer to index the coordinates of the tuples by elements of $\Lambda$.
\begin{proof}
 Assume that the data as given in the statement of the theorem is solvable by a function in $ \mathcal{A}\mathcal{H}(\D^d)$. In view of Theorem 2.1 of \cite{Woerdeman} positive definite matrices $\Gamma_1, \dots, \Gamma_d$ in $\mathcal{B}(\mathbb C^{|\Lambda|})$ exist, with the property:
\begin{align} \label{Eq:Decomposition}
	2(EC^* + CE^*) = \sum_{j=1}^d \Gamma_j -  \sum_{j=1}^d  T_j \Gamma_j T_j^* ,
\end{align} 
where 
\begin{align*}
	C= \operatorname{Col}(C_{\boldsymbol{\beta}})_{\boldsymbol{\beta} \in \Lambda} \ \text{ and } \  E= \operatorname{Col}(\delta_{\boldsymbol{0}   \boldsymbol{\beta} })_{\boldsymbol{\beta} \in \Lambda} 
\end{align*}
are the column operators from $\mathbb{C}$ to $\mathbb{C}^{|\Lambda|}$, 
 $\delta_{\boldsymbol{\beta} \boldsymbol{\gamma} }$ is the Dirac's delta function on $\mathbb{N}_0^d$, and 
\begin{align*}
	T_r = \left(   t^{(r)}_{\boldsymbol{\beta}, \boldsymbol{\gamma}}    \right)_{\boldsymbol{\beta}, \boldsymbol{\gamma} \in \Lambda}
\ \text{ with }
t^{(r)}_{\boldsymbol{\beta}, \boldsymbol{\gamma}}  & = 1, \ \  \text{if } \boldsymbol{\beta}= \boldsymbol{\gamma} + e_{r}; \\
                                                      & = 0, \ \  \text{otherwise}.
\end{align*}
Here $e_r $ is the standard unit vector $(0, \dots, 1, \dots, 0) \in \mathbb{C}^d$ with $1$ at the $r$-th place.
Define, 
\begin{align*}
	X= \left(   c_{\boldsymbol{\beta} - \boldsymbol{\gamma}}   \right)_{  \boldsymbol{\beta} ,  \boldsymbol{\gamma}  \in \Lambda}, \ \text{where we put } c_{\boldsymbol{\beta}}=0 \text{ for }  \boldsymbol{\beta} \notin \Lambda.
\end{align*}

\subsection*{Identification of $\mathbb{C}^{|\Lambda|}$ as a function space}
We interpret $\mathbb{C}^{| \Lambda |}$ as the following subspace of the polynomials in $d$ variables:
\begin{align*}
	\mathbb{C}_{\boldsymbol{n}} [\boldsymbol{z}] = \left \lbrace \sum_{\boldsymbol{\alpha} \leq \boldsymbol{n}} a_{\boldsymbol{\alpha}}  \boldsymbol{z}^{\boldsymbol{\alpha}} : a_{\boldsymbol{\alpha}}  \in \mathbb{C}   \right \rbrace 
\end{align*}
with the usual Euclidean norm $\| \cdot\|_2$.

For every $1\leq j \leq d$, we can view the matrices $T_j$ as the linear operators: 
\begin{align*}
	T_j g (\boldsymbol{z}) = \mathbb{P}_{\boldsymbol{n}} (z_j g) (\boldsymbol{z})  \ \ \text{ for } g\in \mathbb{C}_{\boldsymbol{n}} [\boldsymbol{z}] ,
\end{align*}
where $ \mathbb{P}_{\boldsymbol{n}}$ is the orthogonal projection from the space of hololmorphic polynomials onto the subspace $\mathbb{C}_{\boldsymbol{n}} [\boldsymbol{z}]$.
	In a similar manner the matrix $X$ can be viewed as 
	\begin{align*}
		(Xg)(\boldsymbol{z}) =  \mathbb{P}_{\boldsymbol{n}} (\psi g) (\boldsymbol{z}) \ \ \text{ for } g\in \mathbb{C}_{\boldsymbol{n}} [\boldsymbol{z}] ,
	\end{align*}
	where $\psi$ is the polynomial corresponding to the given interpolation data. Clearly, $X$ is a bounded operator on $\mathbb{C}_{\boldsymbol{n}} [\boldsymbol{z}]$.
	
	\vspace{3mm}
	
	 For any element $h \in \mathbb{C}_{\boldsymbol{n}} [\boldsymbol{z}] $, equation \eqref{Eq:Decomposition} implies
	\begin{align*}
		& 2	\left \langle  CE^*h, h    \right \rangle +  2	\left \langle  EC^*h, h    \right \rangle  = \sum_{j=1}^d \| \Gamma_j^{1/2} h\|^2  - \sum_{j=1}^d \|  \Gamma_j^{1/2} T_j^* h\|^2  \\
			i.e., & \ 4 \Re \left \langle  (X^*h)(0), h(0)  \right \rangle   =  \sum_{j=1}^d \| \Gamma_j^{1/2} h\|^2  - \sum_{j=1}^d \|  \Gamma_j^{1/2} T_j^* h\|^2  \\
			i.e., & \ \| h(0) + (X^*h)(0) \|^2 + \sum_{j=1}^d \|  \Gamma_j^{1/2} T_j^* h\|^2  = \| h(0) - (X^*h)(0) \|^2 + \sum_{j=1}^d \| \Gamma_j^{1/2} h\|^2.
	\end{align*}
		The above identity gives rise to an isometry $W$ defined as 
	\begin{align}\label{Eq:Isometry}
		W \begin{bmatrix}
		  \oplus_{j=1}^d \Gamma_j^{1/2} h \\  h(0) - (X^*h)(0)
		\end{bmatrix}
		= \begin{bmatrix}
			\oplus_{j=1}^d \Gamma_j^{1/2} T_j^*h \\  h(0) + (X^*h)(0)
		\end{bmatrix}
	\end{align}
 acting from the subspace 
\begin{align*}
	\mathcal{M}_1 \bydef \operatorname{span} \left \lbrace   
	 \begin{bmatrix}
		\oplus_{j=1}^d \Gamma_j^{1/2} h \\  h(0) - (X^*h)(0)
	\end{bmatrix} : h \in \mathbb{C}_{\boldsymbol{n}} [\boldsymbol{z}]   \right \rbrace 
	 \subseteq  \left(\oplus_{j=1}^d \mathbb{C}_{\boldsymbol{n}} [\boldsymbol{z}] \right) \oplus \mathbb{C}
\end{align*}
onto the subspace 
\begin{align*}
	\mathcal{M}_2 \bydef \operatorname{span} \left \lbrace   
	\begin{bmatrix}
		\oplus_{j=1}^d \Gamma_j^{1/2} T_j^* h \\  h(0) + (X^*h)(0)
	\end{bmatrix} : h \in \mathbb{C}_{\boldsymbol{n}} [\boldsymbol{z}]   \right \rbrace 
	\subseteq  \left(\oplus_{j=1}^d \mathbb{C}_{\boldsymbol{n}} [\boldsymbol{z}] \right) \oplus \mathbb{C}.
\end{align*}
We extend $W$ to a unitary matrix 
\begin{align*}
	\mathcal{U}= \begin{bmatrix}
		U_{11} & U_{12} \\ U_{21} & U_{22}
	\end{bmatrix} :  \left(\oplus_{j=1}^d \mathbb{C}_{\boldsymbol{n}} [\boldsymbol{z}] \right) \oplus \mathbb{C} \rightarrow  \left(\oplus_{j=1}^d \mathbb{C}_{\boldsymbol{n}} [\boldsymbol{z}] \right) \oplus \mathbb{C}.
\end{align*}
Equation \eqref{Eq:Isometry} yields
\begin{align}  \label{Eq:Unitary}
	\begin{bmatrix}
		U_{11}^* & U_{21}^* \\ U_{12}^* & U_{22}^*
	\end{bmatrix} 
	\begin{bmatrix}
		\oplus_{j=1}^d \Gamma_j^{1/2} T_j^*h \\  h(0) + (X^*h)(0)
	\end{bmatrix}
	=	\begin{bmatrix}
	\oplus_{j=1}^d \Gamma_j^{1/2} h \\  h(0) - (X^*h)(0)
\end{bmatrix}
\end{align}
for every $h$ in $\mathbb{C}_{\boldsymbol{n}} [\boldsymbol{z}]$.

Let us take $\alpha \in \mathbb{C}$ and treat it as a constant polynomial in $\mathbb{C}_{\boldsymbol{n}} [\boldsymbol{z}]$. Note that
\begin{align*}
	T_j^* \alpha =0 \ \  \text{ and } (X^*\alpha)(0)=  \alpha.
\end{align*}
Therefore from \eqref{Eq:Unitary} we have
\begin{align*}
	& U_{22}^* \left( \alpha +    \alpha  \right) =   \alpha  -   \alpha  =0  \implies  U_{22}^* \alpha =0,
\end{align*}
for every $\alpha \in \mathbb{C}$. Thus $U_{22}^* =0$. 

Define $V=-U_{12} \ \ \text{ and } U= U_{11}-U_{12}U_{21}$. Since $\mathcal{U}$ is a unitary, $V$ is an isometry and $U$ is a unitary with $V^*U= U_{21}$. Again, \eqref{Eq:Unitary} implies
\begin{align} \label{Eq:Main}
	\begin{bmatrix}
		U^* & U^*V \\ V^* & 1
	\end{bmatrix}
	\begin{bmatrix}
		\oplus_{j=1}^d \Gamma_j^{1/2} T_j^* h \\ h(0)
	\end{bmatrix} = 
	\begin{bmatrix}
	\oplus_{j=1}^d \Gamma_j^{1/2} h\\ (X^*h)(0)
	\end{bmatrix}.
\end{align}

\noindent \textbf{Claim:} Let \begin{align*}
	\Gamma = [\Gamma_1^{1/2}, \dots, \Gamma_d^{1/2}] : \oplus_{j=1}^d  \mathbb{C}_{\boldsymbol{n}} [\boldsymbol{z}] \rightarrow  \oplus_{j=1}^d  \mathbb{C}_{\boldsymbol{n}} [\boldsymbol{z}] .
\end{align*}
Then 
\begin{align*}
	\Gamma(\xi)(\boldsymbol{z}) = \mathbb{P}_{\boldsymbol{n}} \left( V^*U (I-\Delta(\boldsymbol{z}))U  \right)^{-1} \xi
\end{align*}
for $\xi$ in $\oplus_{j=1}^d  \mathbb{C}_{\boldsymbol{n}} [\boldsymbol{z}]$, where 
$
\Delta(\boldsymbol{z})= \left[ \begin{smallmatrix}
	z_1 I &  & \\ & \ddots &  \\ & & z_d I
\end{smallmatrix} \right]: \oplus_{j=1}^d  \mathbb{C}_{\boldsymbol{n}} [\boldsymbol{z}] \rightarrow  \oplus_{j=1}^d  \mathbb{C}_{\boldsymbol{n}} [\boldsymbol{z}] .
$

\begin{proof}[\textbf{Proof of the claim}] 
	Define, \begin{align*}
		\widetilde{\Gamma} = [\widetilde{\Gamma}_1, \dots, \widetilde{\Gamma}_d] : \oplus_{j=1}^d  \mathbb{C}_{\boldsymbol{n}} [\boldsymbol{z}] \rightarrow  \oplus_{j=1}^d  \mathbb{C}_{\boldsymbol{n}} [\boldsymbol{z}] 
	\end{align*}
	by $\widetilde{\Gamma} (\xi)(\boldsymbol{z}) = \mathbb{P}_{\boldsymbol{n}} \left( V^*U (I-\Delta(\boldsymbol{z})U)^{-1} \right)  \xi$ for $\xi= \oplus_{j=1}^d \xi_j \in \oplus_{j=1}^d  \mathbb{C}_{\boldsymbol{n}} [\boldsymbol{z}] $.
	
Consequently, \begin{align} \label{Eq:Mainak}
		& [T_1 \widetilde{\Gamma}_1, \dots, T_d \widetilde{\Gamma}_d] \left( \oplus_{j=1}^d \xi_j \right)(\boldsymbol{z})  \notag \\
		& = \sum_{j=1}^d T_j \widetilde{\Gamma}_j (0 \oplus \dots \oplus \underbrace{\xi_j}_{\text{j-th position}} \oplus \dots \oplus 0)(\boldsymbol{z}) \notag \\
		&=  \mathbb{P}_{\boldsymbol{n}} \left(    \sum_{j=1}^d z_j \widetilde{\Gamma}_j (0 \oplus \dots \oplus \underbrace{\xi_j}_{\text{j-th position}} \oplus \dots \oplus 0)(\boldsymbol{z})    \right)  \notag  \\
		& = \mathbb{P}_{\boldsymbol{n}} \left(   [\widetilde{\Gamma}_1, \dots, \widetilde{\Gamma}_d ] \Delta(\boldsymbol{z}) (\oplus_{j=1}^d \xi_j)(\boldsymbol{z})  \right)  \notag \\
		&= \mathbb{P}_{\boldsymbol{n}} \left(   V^*U (I- \Delta(\boldsymbol{z})U)^{-1} \Delta(\boldsymbol{z})  \right ) (\oplus_{j=1}^d \xi_j) .
	\end{align}
	
	Now, for $h \in \mathbb{C}_{\boldsymbol{n}} [\boldsymbol{z}] $ and $\xi \in \oplus_{j=1}^d \mathbb{C}_{\boldsymbol{n}} [\boldsymbol{z}] $,
	\begin{align*}
	&	\langle  h, \widetilde{\Gamma} \xi \rangle \\
	=& \left \langle  h, \mathbb{P}_{\boldsymbol{n}} \left(   V^*U (I- \Delta(\boldsymbol{z})U)^{-1} (I- \Delta(\boldsymbol{z}) U) \xi  \right )   \right \rangle
	+\left \langle  h, \mathbb{P}_{\boldsymbol{n}} \left(   V^*U (I- \Delta(\boldsymbol{z})U)^{-1} \Delta(\boldsymbol{z}) U \xi  \right )   \right \rangle.
	\end{align*}
	Again \eqref{Eq:Mainak} yields:
	\begin{align}
		& \left \langle  h, \widetilde{\Gamma} \xi \right \rangle    
		=  \left \langle  h,  \mathbb{P}_{\boldsymbol{n}} (V^*U\xi) \right \rangle_{\mathbb{C}_{\boldsymbol{n}} [\boldsymbol{z}]}   + \left \langle  h,  [T_1 \widetilde{\Gamma}_1, \dots, T_d \widetilde{\Gamma}_d ] U\xi \right \rangle_{\mathbb{C}_{\boldsymbol{n}} [\boldsymbol{z}]}   \notag  \\
		=& \left \langle  U^*Vh(0), \xi \right \rangle + \left \langle    U^*(\oplus_{j=1}^d \widetilde{\Gamma}_j^* T_j^* h), \xi    \right \rangle.
	\end{align}
	
	Therefore, \begin{align} \label{Eq:Mihai}
		\oplus_{j=1}^d \widetilde{\Gamma}_j^* h = U^*Vh(0) + U^* \left( \oplus_{j=1}^d \widetilde{\Gamma}_j^* T_j^* h  \right).
	\end{align}
	Combining \eqref{Eq:Main} and \eqref{Eq:Mihai}, we obtain
	\begin{align*}
		\oplus_{j=1}^d \left( \Gamma_j^{1/2}  - \widetilde{\Gamma}_j \right)h = U^* \oplus_{j=1}^d    \left( \Gamma_j^{1/2}  - \widetilde{\Gamma}_j^* \right)T_j^*h.
	\end{align*}
	
	We use Lemma 2.2 of \cite{Ball-Li-Timotin-Trent} to get
	\begin{align*}
		\widetilde{\Gamma}_j^* = \Gamma_j^{1/2} \ \ \text{ for } 1\leq j \leq d.
	\end{align*}
	Consequently, for $\xi = \oplus_{j=1}^d \xi_j $,
$$
		[\Gamma_1^{1/2}, \dots, \Gamma_d^{1/2}](\xi) (\boldsymbol{z}) = \mathbb{P}_{\boldsymbol{n}} \left(  V^*U(I - \Delta(\boldsymbol{z})U)^{-1} \xi \right)
	$$
	and 
	\begin{align}\label{Eq:Mainak-1}
		[T_1 \Gamma_1^{1/2}, \dots, T_d \Gamma_d^{1/2}] (\xi)(\boldsymbol{z}) =  \mathbb{P}_{\boldsymbol{n}} \left(  V^*U(I - \Delta(\boldsymbol{z})U)^{-1} \Delta(\boldsymbol{z}) \xi \right).
	\end{align}
	This proves our claim.
\end{proof}
	\subsection*{Concluding part} We return to the proof of Theorem \ref{Thm:Inner-solution}.
	Consider 
	\begin{align*}
		\varphi(\boldsymbol{z}) = 1 + 2 V^* U (I - \Delta(\boldsymbol{z})U)^{-1} \Delta(\boldsymbol{z}) V.
	\end{align*}
	The above representation of $\varphi$ is known as a realization formula for $\varphi$. In particular, $\varphi \in \mathcal{A} \mathcal{H}(\D^d)$. 
For $h \in \mathbb{C}_{\boldsymbol{n}} [\boldsymbol{z}] $ and $\alpha \in \mathbb{C}$,
	\begin{align*}
		& \langle  h, T_{\varphi} \alpha  \rangle \\
		=& \langle  h, \alpha + 2V^*U(I- \Delta(\boldsymbol{z})U)^{-1} \Delta(\boldsymbol{z}) V\alpha \rangle \\
		=& \langle h(0), \alpha \rangle + 2 \langle  Vh, U(I- \Delta(\boldsymbol{z})U)^{-1} \Delta(\boldsymbol{z}) V\alpha \rangle \\
		=&   \langle h(0), \alpha \rangle + 2 \langle   V^*\left(\oplus_{j=1}^d \Gamma_j^{1/2} T_j^* h\right), \alpha  \rangle \ \ (\text{by }\eqref{Eq:Mainak-1}) \\
		=&   \langle h(0), \alpha \rangle + \langle   (X^*h)(0) - h(0) , \alpha \rangle  \ \ (\text{using } \eqref{Eq:Main})  \\
		=& \langle  (X^*h)(0), \alpha \rangle = \langle h, X\alpha \rangle.  
	\end{align*}

\noindent Suppose, $
\varphi(\boldsymbol{z})= \sum_{\boldsymbol{\beta} \in \mathbb{N}_0^d} \varphi_{\boldsymbol{\beta}} \boldsymbol{z}^{\boldsymbol{\beta}}
$ is the Taylor series of $\varphi$ at the origin.
Then for $\boldsymbol{\beta} \in \Lambda$,
\begin{align*}
	 \left \langle   T_\varphi^* h ,  M_{\boldsymbol{z}}^{\boldsymbol{\beta}}  \alpha \right \rangle = \left \langle   
	 T_\varphi^* M_{\boldsymbol{z}}^{* \boldsymbol{\beta}}  h ,   \alpha \right \rangle 
	 = \left \langle   X^* T^{* \boldsymbol{\beta}}  h ,   \alpha \right \rangle 
	 = \left \langle  X^* h ,   M_{\boldsymbol{z}}^{\boldsymbol{\beta}}  \alpha  \right \rangle.
\end{align*}
	Therefore, we have $\varphi_{\boldsymbol{\beta}}= c_{\boldsymbol{\beta}}$ for every $\boldsymbol{\beta} \in \Lambda$.
Thus $\varphi $ is a solution of the Carath\'eodory-Fej\'er problem with the given interpolation data. To see that, $\varphi$ is a rational function with poles off $\D^d$, we look at the expression 
\begin{align*}
	\varphi(\boldsymbol{z}) = 1 + 2 V^*U (I- \Delta(\boldsymbol{z})U)^{-1} \Delta(\boldsymbol{z}) V \ \ \text{ for } \boldsymbol{z} \in \D^d,
\end{align*}
where $U$  is a unitary matrix of order $d |\Lambda|$ and $V$ is a row-matrix of size $d |\Lambda| \times 1$. Therefore,
\begin{align*}
	\varphi(\boldsymbol{z}) = 1 +  \frac{1}{\operatorname{det}(I- \Delta(\boldsymbol{z})U)}  2V^*U  \operatorname{Adj}(I- \Delta(\boldsymbol{z})U)  \Delta(\boldsymbol{z}) V
\end{align*}
where $\operatorname{det}(I- \Delta(\boldsymbol{z})U)$ is non-zero in $\D^d$ as $\Delta(\boldsymbol{z})U$ is a strict contraction. This proves that $\varphi$ is a rational function with poles off $\D^d$. 

Moreover, we can show that $\Re \varphi(\boldsymbol{z})$ is zero $\Theta$-a.e. on $\mathbb{T}^d$. Indeed, for $\boldsymbol{z} \in \D^d$,
\begin{align} \label{Eq:Cayley-inner}
&2 \Re \varphi(\boldsymbol{z})  \notag \\
=&2 V^* \left[ 1+ U (I- \Delta(\boldsymbol{z})U)^{-1}  \Delta (\boldsymbol{z}) +  \Delta (\boldsymbol{z}) ^*  (I- U^* \Delta(\boldsymbol{z})^*)^{-1}U^*  \right] V   \notag\\
=& 2 V^*   \left[ 1+ U (I- \Delta(\boldsymbol{z})U)^{-1}  U \Delta (\boldsymbol{z}) +  \Delta (\boldsymbol{z}) ^* U^*   (I- \Delta(\boldsymbol{z})^* U^* )^{-1} \right] V   \notag\\
=& 2 V^* (I- \Delta(\boldsymbol{z})U)^{-1}  \left[    I- U \Delta(\boldsymbol{z})   \Delta(\boldsymbol{z}) ^* U^*        \right] (I- \Delta(\boldsymbol{z})^* U^* )^{-1} V.
\end{align}

The rational function $\varphi$ can have poles on $\mathbb{T}^d$ only in a measure zero set which is contained in the zero set of the polynomial $\operatorname{det}(I - \Delta(\boldsymbol{z})U)$. Therefore, for almost every point $\boldsymbol{\zeta}$ in $\mathbb{T}^d$, the radial limit 
\begin{align*}
	& 2 \Re \varphi(r \boldsymbol{\zeta}) \\
	=& 2V^* (I- r \Delta(\boldsymbol{\zeta})U)^{-1}  \left[    I- r^2 U \Delta(\boldsymbol{\zeta})   \Delta(\boldsymbol{\zeta}) ^* U^*        \right] (I- r \Delta(\boldsymbol{\zeta})^* U^* )^{-1} V  \\
	=& 2V^* (I- r \Delta(\boldsymbol{\zeta})U)^{-1}  \left[    I- r^2   \right] (I- r \Delta(\boldsymbol{\zeta})^* U^* )^{-1} V  \rightarrow 0 
\end{align*}
as $r \rightarrow 1^{-}$. In conclusion, $\varphi$  is a Cayley rational inner solution of the given (solvable) Carath\'eodory-Fej\'er interpolation data.
\end{proof}

\section{The coefficient sets} \label{The coefficient problem}

We continue to work in the polydisk $\D^d$. For a multi-index $\boldsymbol{n} $ in $\mathbb{N}_0^d$, we consider the set of normalized Taylor coefficients of a function belonging to the Agler's class, $ \mathcal{A}\mathcal{H}(\D^d)$ as follows:
$$ K_{\boldsymbol{n}} = \{ (c_{\boldsymbol{\alpha}}(\varphi))_{\boldsymbol{\alpha} \in \Lambda_{\boldsymbol{n}}} :  \varphi \in  \mathcal{A}\mathcal{H}(\D^d) \text{ with } c_{\boldsymbol{0}}(\varphi) = 1\}.$$
We owe to Pfister a first account of some basic properties of this set of coefficients for the Agler class functions in the Schur class setting \cite{Pfister}. From his analysis, we observe that for every fixed $\boldsymbol{n}$, $K_{\boldsymbol{n}}$ is a closed convex set, with coefficients associated to 
a Cayley rational function belonging to the boundary of $K_{\boldsymbol{n}}.$  A full description of $K_{\boldsymbol{n}}$ seems to be out of reach for $d \geq 2$. However, Woerdeman's matrix realization
approach to Agler's class of functions implies the following observation. Since we shall use some results from the theory of real algebraic geometry, we view the set $K_{\boldsymbol{n}}$ as a subset of $\R^{2  |\Lambda_{\boldsymbol{n}}|}$ by splitting the elements into their real and imaginary parts.
\begin{theorem}\label{semialgebraic}
For any $d \geq 1$ and $\boldsymbol{n} \in \mathbb{N}_0^d$ the set $K_{\boldsymbol{n}} \subset \R^{2  |\Lambda_{\boldsymbol{n}}|}$ is semi-algebraic.
\end{theorem}
\begin{proof} Recall that, a subset of the Euclidean space $\R^n$ is called semi-algebraic if it is a finite union of sets defined by real polynomial equalities and inequalities, cf. \cite[Chapter 2]{BCR}.
We interpret Woerdeman's matricial identity (\ref{Eq:Decomposition}) as a certificate of membership $(c_{\boldsymbol{\alpha}})_{ \boldsymbol{\alpha} \in \Lambda_{\boldsymbol{n}}} \in K_{\boldsymbol{n}}$  expressed in the form of a system of real polynomial equations with additional variables. For instance, the positivity of a finite matrix $\Gamma$ can be encoded into the representation $\Gamma = X^\ast X$, and so on.

In view of Tarski-Seidenberg Principle stated in its geometric form (a projection of a semi-algebraic set is semi-algebraic) we find that $K_{\boldsymbol{n}}$ is semi-algebraic. See Theorem 2.2.1 in \cite{BCR} and the adjacent discussion of stability of semi-algebraic sets under algebraic or topological operations. In other terms, all additional variables appearing on the real system of equations encoded by (\ref{Eq:Decomposition}) can be eliminated, leading to a characterization of the coefficient collection  $K_{\boldsymbol{n}}$ as a finite union of sets of equalities and inequalities involving solely the variables $(c_{\boldsymbol{\alpha}})_{ \boldsymbol{\alpha} \in \Lambda_{\boldsymbol{n}}}$.
\end{proof}

The {\em algebraic boundary} of a semi-algebraic set is the smallest algebraic variety containing its boundary in the euclidean topology. For further details on real algebraic geometry we refer to \cite{Sinn}.
\begin{corollary} For any multi-index $\boldsymbol{n} $ in $\mathbb{N}_0^d$, the algebraic boundary of the coefficient set  $K_{\boldsymbol{n}}$ is a hypersurface.
\end{corollary}

\begin{proof} We have already observed that the set  $K_{\boldsymbol{n}}$ is closed and convex. We show below that $K_{\boldsymbol{n}}$ has a non-empty interior in the Euclidean topology of  $\R^{2  |\Lambda_{\boldsymbol{n}}|}$. 
	
	Consider the linear map 
	$$
			T: \mathcal{F}_d \rightarrow \{(c_{\boldsymbol{\alpha}}(f))_{\boldsymbol{\alpha} \in \Lambda_{\boldsymbol{n}}}: f \in \mathcal{F}_d \} \subseteq \mathbb{C}^{|\Lambda_{\boldsymbol{n}}|} \  \text{ defined as} \\
	$$
	$$
	T(f) =  (c_{\boldsymbol{\alpha}}(f))_{\boldsymbol{\alpha} \in \Lambda_{\boldsymbol{n}}} \ \text{ for } f\in \mathcal{F}_d.
	$$
For every $f\in \mathcal{F}_d$, $\| f\|_{\operatorname{Agler}} \leq \|f\|_\infty$ and so by the Cauchy's estimates we see that $T$ is continuous. 
Note that for any standard unit vector $e_{\boldsymbol{\alpha}} =(0, \dots, 1, \dots, 0) \in \mathbb{C}^{|\Lambda_{\boldsymbol{n}}|}$, we take $\boldsymbol{z}^{\boldsymbol{\alpha}}$ which is obviously in $\mathcal{F}_d$ and $T(\boldsymbol{z}^{\boldsymbol{\alpha}}) = e_{\boldsymbol{\alpha}} $. Thus, $T$ being linear it is surjective and hence an open map. Therefore, $T$ maps the open Agler-norm unit ball to an open set, $U$ containing the null vector in $\mathbb{C}^{|\Lambda_{\boldsymbol{n}}|}$.

\noindent Now, recall the universal rational functions $M_{\boldsymbol{\alpha}}$ described in Subsection \ref{Subsection: Change of Taylor coeff.}. We consider the rational map 
$$(M_{\boldsymbol{\alpha}})_{\boldsymbol{\alpha} \leq \boldsymbol{n}} : \{(c_{\boldsymbol{\alpha}}(\varphi))_{\boldsymbol{\alpha} \in \Lambda_{\boldsymbol{n}}}: \varphi \in \mathcal{A}\mathcal{H}(\D^d) \}  \rightarrow \{(c_{\boldsymbol{\alpha}}(f))_{\boldsymbol{\alpha} \in \Lambda_{\boldsymbol{n}}}: f \in \mathcal{A}\mathcal{S}(\D^d) \} .
$$
Then this map is continuous as it has no singularity as $(1+ c_{\boldsymbol{0}}(\varphi)) \neq 0$ for every $\varphi$ in $\mathcal{A}\mathcal{H}(\D^d)$. Therefore, the inverse image, $\tilde{U}$ of the open set $U$ under this rational map is open. Also, the inverse image of the intersection of $U$ and the hyperplane $\{(c_{\boldsymbol{\alpha}}(f))_{\boldsymbol{\alpha} \in \Lambda_{\boldsymbol{n}}}: f \in \mathcal{A}\mathcal{S}(\D^d) \text{ and } f(0)=0 \}$ is the same as $\tilde{U} \cap K_{\boldsymbol{n}}$ which contains $(c_{\boldsymbol{\alpha}}(\mathbf{1}))_{\boldsymbol{\alpha} \in \Lambda_{\boldsymbol{n}}}$. Here $\mathbf{1}$ denotes the constant one function.
 Thus $K_{\boldsymbol{n}}$ contains a non-empty open set in the subspace topology of $\mathbb{C}^{|\Lambda_{\boldsymbol{n}}|}$ which is the same as that of $\mathbb{R}^{2|\Lambda_{\boldsymbol{n}}|}$.
 
 For any $\varphi \in K_{\boldsymbol{n}}$, from the Kor\'anyi-Puk\'anszky integral representation \eqref{Eq:KP-rep} we have 
 $$
 \varphi(\boldsymbol{z}) = \int_{\mathbb{T}^d}  \left(  \frac{2}{\prod_{j=1}^{d} (1 - z_j \bar{\xi}_j)}-1  \right) d\nu(\boldsymbol{\xi}),
 $$
 for a unique regular Borel probability measure $\nu$ on $\mathbb{T}^d$ satisfying the vanishing moment conditions \eqref{Eq:Pluri-measure}. Therefore, 
 $$
 c_{\boldsymbol{0}}(\varphi) = 1  \ \text{ and } \ c_{\boldsymbol{\alpha}}(\varphi) = 2 \widehat{\nu}(\boldsymbol{\alpha}) \text{ for } \boldsymbol{\alpha} \in \mathbb{N}_0^d \setminus \{ \boldsymbol{0}\}.
 $$
 Since the total mass of the measure $\nu$ is $1$, the above expressions yields that $\|c_{\boldsymbol{\alpha}}(\varphi) \| \leq 2$ for every $\boldsymbol{\alpha} \in  \mathbb{N}_0^d$. Thus $K_{\boldsymbol{n}}$ is compact. 
 Therefore, Corollary 2.8 of \cite{Sinn} implies that the algebraic boundary of $K_{\boldsymbol{n}}$ is a hypersurface.
\end{proof}

Notice that the coefficient sets are nested: if $\boldsymbol{n} \leq \boldsymbol{m}$ entrywise, then the projection of $K_ {\boldsymbol{m}} \subset \R^{2  |\Lambda_{\boldsymbol{m}}|}$ onto $\R^{2  |\Lambda_{\boldsymbol{n}}|}$
is $K_{\boldsymbol{n}}$. The description of the coefficients sets in 1D is classical, involving determinants of Toeplitz or Hankel matrices cf. \cite{Ahiezer-Krein}.

\subsection*{Invariance of the coefficient sets }
The full set of admissible coefficients
$$
\widetilde{K}_{\boldsymbol{n}} = \{ (c_{\boldsymbol{\alpha}}(\varphi))_{ \boldsymbol{\alpha} \in \Lambda_{\boldsymbol{n}}} :  \varphi \in  \mathcal{A}\mathcal{H}(\D^d)\}
$$
is the convex cone generated by  $K_{\boldsymbol{n}}$.

Let $\boldsymbol{c}= (c_{\boldsymbol{\alpha}})_{\alpha \in \Lambda_{\boldsymbol{n}}}$ be a point in $\widetilde{K}_{\boldsymbol{n}}$ so that there exists $\varphi \in \mathcal{A}\mathcal{H}(\D^d)$ such that $ c_{\boldsymbol{\alpha}}= c_{\boldsymbol{\alpha}}(\varphi) $ for every $\boldsymbol{\alpha} \in \Lambda_{\boldsymbol{n}}$. 
Suppose $\psi$ is an automorphism of the right half plane (i.e., $\psi$ is a Mobius transformation taking the right half plane onto itself). Then $\psi \circ \varphi$ is in $\mathcal{A}\mathcal{H}(\D^d)$. Indeed, for a $d$-tuple of strict contractions, $\underline{T}= (T_1, \dots, T_d)$ on a fixed infinite dimensional Hilbert space, $\Re \varphi(\underline{T}) \geq 0$ implies that $\| \eta (\varphi(\underline{T})) \| \leq 1$, where $\eta(z) = (z-1)/(z+1)$ for $z$ in the right half plane. Now, $\eta \circ \psi \circ \eta^{-1}: \D \rightarrow \D$ is holomorphic in $\D$ and hence by the von Neumann inequality,
$$
\| \eta \circ \psi \circ \eta^{-1} \circ( \eta \circ \varphi(\underline{T})) \| \leq 1
$$
that is, $\eta \circ \psi (\varphi(\underline{T}))$ is a contraction. Thus $\Re \psi (\varphi(\underline{T}) \geq 0$. 

By an argument as in Subsection  \ref{Subsection: Change of Taylor coeff.} we observe that the Taylor coefficient $c_{\boldsymbol{\alpha}}( \psi \circ \varphi)$ is a rational function $ L_{\boldsymbol{\alpha}}$ (say), depending the derivatives of $\psi$ at $0$ up to order $|\boldsymbol{\alpha}|$, with the Taylor coefficients $c_{\boldsymbol{\beta}}(\varphi)$ with $\boldsymbol{\beta} \leq \boldsymbol{\alpha}$ as the variables.
Therefore, $(L_{\boldsymbol{\alpha}})_{\boldsymbol{\alpha} \in \Lambda_{\boldsymbol{n}}}$ is in $\widetilde{K}_{\boldsymbol{n}}$.

Thus for every $\psi$ and $\boldsymbol{c}$ as above, we write $L_{\boldsymbol{\alpha}} = L_{\boldsymbol{\alpha}} \left((c_{\boldsymbol{\beta}})_{\boldsymbol{\beta} \leq \boldsymbol{\alpha}}, \psi \right)$ to emphasize the dependence on $\psi$ and $\boldsymbol{c}$.  The fact that $\left( L_{\boldsymbol{\alpha}} \right)_{\boldsymbol{\alpha} \in \Lambda_{\boldsymbol{n}}} \in \widetilde{K}_{\boldsymbol{n}}$ reflects the {\em invariance of $\widetilde{K}_{\boldsymbol{n}}$ under the automorphism $\psi$}.

In particular, if $\boldsymbol{c} \in K_{\boldsymbol{n}}$ and $\psi(1)=1$, then $\left( Q_{\boldsymbol{\alpha}} \right)_{\boldsymbol{\alpha} \in \Lambda_{\boldsymbol{n}}} \in K_{\boldsymbol{n}}$.

\subsection*{Example}
For any automorphism $\Phi$ of the right half plane with $\Phi(1)=1$, a point 
$( 1, c_{01}, c_{10}, c_{11}) $ is in $K_{11} $ if and only if 

\noindent $(1, \Phi'(c_{00}) 	c_{01}, \Phi'(c_{00}) 	c_{10}, \Phi''(c_{00}) c_{01} c_{10}  + \Phi'(c_{00}) c_{11})$ is in $K_{11}$.

\subsection*{Description of $K_{11}$ in the case $\D^2$}

As a verification of Theorem \ref{semialgebraic} we show that in the simplest low degree/low dimension case, accessible by direct computation, the set of admissible coefficients is semi-algebraic.
Pfister's doctoral dissertation \cite{Pfister-thesis} contains a description of the coefficient set $K_{11}$ in the two variable Schur class setting. We sketch below a proof of the companion result 
for the Herglotz-Nevanlinna class. 

Let $(1, c_{01}, c_{10}, c_{11})$ be a point in $K_{11}$.  Since the Carth\'eodory-Fej\'er interpolation problem with data $\{ 1, c_{01}, c_{10}, c_{11}\}$ has a solution in the Herglotz-Nevanlinna class, the matrix $X+X^*$ must be non-negative (see \cite[Theorem 2.1]{Woerdeman}), where 
\begin{align*}
	X= \begin{bmatrix}
		1 & 0 & 0 & 0 \\ c_{01} & 1 & 0 & 0 \\ c_{10} & 0 & 1 & 0 \\ c_{11} & c_{10} & c_{01} & 1 
	\end{bmatrix}.
\end{align*}
Therefore, \begin{align*}
	X + X^*= \begin{bmatrix}
		2 & \overline{c_{01}} & \overline{c_{10}} &  \overline{c_{11}} \\ c_{01} & 2 & 0 & \overline{c_{10}} \\ c_{10} & 0 & 2 & \overline{c_{01}} \\ c_{11} & c_{10} & c_{01} & 2
	\end{bmatrix} \geq 0
\end{align*}
and so,
\begin{align} \label{Eq:cond-1}
	2|c_{11} - c_{10} c_{01} | + |c_{10}|^2 + |c_{01}|^2 \leq 4.
\end{align}

Let $f$ be a Herglotz-Nevanlinna function on $\D^2$ such that 
\begin{align*}
	f(z, w)= 1+ c_{01} w + c_{10} z + c_{11} zw + R(z,w) \ \text{ for } (z,w) \in \D^2, 
\end{align*}
where $R(z,w)$ is the remainder of the Taylor series of $f$ around the origin. For any $\lambda \in \mathbb{T}$, 
\begin{align*}
	f_\lambda (z) \bydef f(z, \lambda z) = 1+ (\lambda c_{01} + c_{10}) z + o(z^2)
\end{align*}
is a Herglotz-Nevanlinna function on $\D$. So, $f_\lambda$ solves the Carath\'eodory-Fej\'er interpolation problem in $\D$ with data $\{1, c_{10}+ \lambda c_{01}\}$. Therefore, the corresponding matrix
$$
\begin{bmatrix}
	2 & \overline{c_{10}}+ \bar{\lambda} \overline{ c_{01}} \\ c_{10}+ \lambda c_{01} & 2 
\end{bmatrix}
$$
is non-negative and hence, $|c_{10}+ \lambda c_{01}| \leq 2.$ Since this inequality is true for every $\lambda \in \mathbb{T}$, we must have
\begin{align} \label{Eq:cond-2}
	|c_{10}|+  |c_{01}| \leq 2.
\end{align}

\noindent Thus the realtions \eqref{Eq:cond-1} and \eqref{Eq:cond-2} are necessary for $( 1, c_{01}, c_{10}, c_{11})$ to be in $K_{11}$.

\begin{proposition} \label{normalized-K_11}
	The point $( 1, c_{01}, c_{10}, c_{11})$ belongs to $K_{11}$ if and only if \\
	$2|c_{11} - c_{10} c_{01} | + |c_{10}|^2 + |c_{01}|^2 \leq 4$ and $|c_{10}|+  |c_{01}| \leq 2.$
	\end{proposition}
	
	\begin{proof}
	We have shown the necessity of the conditions \eqref{Eq:cond-1} and \eqref{Eq:cond-2} just before the statement of the theorem. To prove the sufficiency, first we deal with the case $|c_{01}|^2 + |c_{10}|^2 < 4 $.
	Take $\sigma = \frac{2(c_{11}-c_{10} c_{01})}{4-|c_{10}|^2 -|c_{01}|^2 }$ and 
	$$
	\varphi(z,w) = \frac{ 1+ \frac{\sigma \overline{ c_{01}} + c_{10}}{2} z +   \frac{\sigma \overline{ c_{10}} + c_{01}}{2} w + \sigma zw }{ 1+ \frac{\sigma \overline{ c_{01}} - c_{10}}{2} z +   \frac{\sigma \overline{ c_{10}} - c_{01}}{2} w - \sigma zw  }. 
	$$
The rational function $\varphi$ can be rewritten as 
$$
\varphi(z,w) = \frac{1+ g(z,w)}{1-g(z,w)}
$$
where 
$$
g(z,w)= \frac{\frac{c_{10}}{2} z + \frac{c_{01}}{2} w + \sigma zw}{  1+ \sigma \left(  \frac{\overline{ c_{01}}}{2} z +   \frac{\overline{ c_{10}}}{2} w  \right)  }.
$$
For $(z,w) \in \D^2$, using \eqref{Eq:cond-1} and \eqref{Eq:cond-2} we observe that the denominator of $g$ is non-vanishing. Thus $g$ defines a holomorphic function in $\D^2$. For $0<r <1$, the function $\varphi_r(z,w) \bydef \varphi(rz, rw)$ is continuous in a neighbourhood of $\overline{\D}^2$. Moreover, when $|z|=|w|=1$,
\begin{align*}
	|g_r(z, w)| = |g(rz, rw)| = r  \left| \frac{\frac{c_{10}}{2} z + \frac{c_{01}}{2} w + r \sigma zw}{  1+  r \sigma zw \left(  \frac{\overline{ c_{01}}}{2} \bar{z} +   \frac{\overline{ c_{10}}}{2} \bar{w}  \right)  } \right| \leq r.
\end{align*}
Hence by maximum modulus principle, $|g_r| \leq r $ and thus $|g| \leq 1$ on $\D^2$ yielding $\Re \varphi \geq 0$ in $\D^2$.

On the other hand, if $|c_{01}|^2 + |c_{10}|^2 = 4 $ then \eqref{Eq:cond-2} implies that $c_{01} c_{10} = 0$. Without loss of generality, we assume that $c_{10}=0$. Then $c_{01}=2 \tau $, for some unimodular constant $\tau$. Also, \eqref{Eq:cond-1} implies that $c_{11}=0$. Let us take $\psi(z,w) = \frac{1+\tau w}{1- \tau w}$ on $\D^2$ and clearly $\psi$ is in the Herglotz-Nevanlinna class solving the interpolation data $\{ 1,2\tau, 0, 0\}$. 
This completes the proof of the theorem.
\end{proof}

\begin{corollary}
	The Carath\'eodory-Fej\'er interpolation data $( c_{00}, c_{01}, c_{10}, c_{11})$ is solvable if and only if we the following two relations are satisfied:
	\begin{align}
		&	2| \left( \Phi''(c_{00}) - ( \Phi'(c_{00}) )^2  \right)c_{01} c_{10}  + \Phi'(c_{00}) c_{11}  | + |\Phi'(c_{00})|^2 \left( |c_{01}|^2 + |c_{10}|^2 \right)  \leq 4    \label{Eq:char-1}\\
		& \text{ and }  \ \	|\Phi'(c_{00})|  \left( |c_{01}| + |c_{10}|  \right)  \leq 2 , \label{Eq:char-2}
	\end{align}
	where $\Phi$ is an automorphism of the right half plane given by 
	\begin{align*}
		\Phi(z) = \frac{Az + B }{ Cz + D},
	\end{align*}
	with $A= 2-\gamma - \bar{\gamma}, B=  \bar{\gamma} - \gamma, C=(\gamma -\bar{\gamma}), D= 2 + \gamma + \bar{\gamma}$ and $\gamma = \frac{c_{00}-1}{c_{00} +1}$.
\end{corollary}
\begin{proof}
	Suppose there exists a Herglotz-Nevanlinna function $f$ on $\D^2$ such that $f$ solves the Carath\'eodory-Fej\'er interpolation data $( c_{00}, c_{01}, c_{10}, c_{11})$. Let $ \gamma = \frac{c_{00}-1}{c_{00} +1}$. Consider the automorphism $\Phi $ of the right half plane given by 
	$$
		\Phi(z) = \frac{Az + B }{ Cz + D},
	$$
	where $A= 2-\gamma - \bar{\gamma}, B=  \bar{\gamma} - \gamma, C=(\gamma -\bar{\gamma})$ and $D= 2 + \gamma + \bar{\gamma}$.
	Then $\Phi(c_{00}) =1 $. Note that, the function $g = \Phi \circ f$ is in the Herglotz-Nevanlinna class of $\D^2$ with $g_{00} = g(0,0)=1$. Now,
	\begin{align*}
		g_{01}= \frac{\partial g}{\partial w} (0, 0)= \Phi'(c_{00}) 	c_{01}, \ \ 	g_{10} = \frac{\partial g}{\partial z} = \Phi'(c_{00}) 	c_{10} , \ \text{ and }  \\
		g_{11}= \frac{\partial^2 g}{\partial w \partial z}(0,0)= \Phi''(c_{00}) c_{01} c_{10}  + \Phi'(c_{00}) c_{11}. 
	\end{align*}
	Therefore, $(1, g_{01}, g_{10}, g_{11} ) \in K_{11}$. Applying Theorem \ref{normalized-K_11} to this point in $K_{11}$, we obtain the relations \eqref{Eq:char-1} and \eqref{Eq:char-2}. Remark that the expressions 
	$$\Phi'(c_{00})= \frac{AD-BC}{(Cc_{00}+D)^2}, \ \ \ \Phi''(c_{00})= -\frac{2(AD-BC)C}{(Cc_{00} + d)^3 }$$ are rational functions in $\Re c_{00}$ and $\Im c_{00}$ and hence the inequalities  stated in 
	the corollary yield two polynomial inequalities. 
	
	Conversely, suppose the relations \eqref{Eq:char-1} and \eqref{Eq:char-2} hold for $(c_{00}, c_{01}, c_{10}, c_{11})$.  Therefore by Theorem \ref{normalized-K_11},  $(1, \Phi'(c_{00}) c_{01}, \Phi'(c_{00}) 	c_{10}, \Phi''(c_{00}) c_{01} c_{10}  + \Phi'(c_{00}) c_{11})$ is in $K_{11}$ where $\Phi$ is as in the previous part of the proof. Now we use the invariance of $K_{11}$ under the automorphism $\Phi^{-1}$ on the right half plane to see that $(c_{00}, c_{01}, c_{10}, c_{11}) \in K_{11}$. The proof is complete.
\end{proof}

\section{Rational approximation} \label{Rational approximation}

This section focuses on rational approximation aspects of Herglotz-Nevanlinna, or bounded holomorphic functions defined in the polydisk.

\subsection{Rational approximation in Agler's class}
The following result, of a rather theoretical value only,  is a direct consequence of Theorem \ref{Thm:Inner-solution}. 
\begin{theorem}\label{CF-approx}
	Let $d \geq 1$ and $f: \D^d \rightarrow \overline{\mathbb{D}}$ be in the Agler class, $\mathcal{A}\mathcal{S}(\D^d)$. Then there exists a sequence of rational inner functions $\{ f_n\}_{n\in \mathbb{N}}$ such that:
	\begin{itemize}
	\item{the multi-degree of $f_n$ is at most $(2(n+1)^d, \dots, 2(n+1)^d)$;}
	
	\item{Taylor's series at $\boldsymbol{z}=0$ of $f$ and $f_n$ are the same up to the multi-degree $(n, \dots, n)$ for every $n\in \mathbb{N}$;}
	
         \item{$f_n$ converges to $f$ uniformly on compact subsets of $\D^d$.}
         \end{itemize}
\end{theorem}

To simplify terminology, the uniform convergence on compact subsets of $\D^d$ will be recognized as convergence with respect to the Fr\'echet topology of the space
${\mathcal O}(\D^d)$ of holomorphic functions on $\D^d$. It is important to recall that when $d=2$, $\mathcal{A}\mathcal{S}(\D^d) = \mathcal{S}(\D^d)$. 
\begin{proof}
	Let $f$ be in $\mathcal{A}\mathcal{S}(\D^d)$ as above. If $f$ is identically equal to $1$ then there is nothing to prove. So, we assume that $f$ is not identically $1$. Then $g \bydef (1+f)(1-f)^{-1}$ is in $\mathcal{A}\mathcal{H}(\D^d)$. We take $\boldsymbol{n}=(n,\dots, n)$ for each $n \in \mathbb{N}$. 
	Then by Theorem \ref{Thm:Inner-solution}, for each $\boldsymbol{n}$ we find a Cayley rational inner function $g_n$ such that the Taylor polynomials of $g$ and $g_n$ are the same up to order $\boldsymbol{n}$. 
	 Define, $f_n = (g_n -1)(1+ g_n)^{-1}$. Then the functions $f_n$ are rational inner. From Subsection \ref{Subsection: Change of Taylor coeff.}, we get the rational functions $ Q_{\boldsymbol{\alpha}} $ for $ \boldsymbol{\alpha} \in \mathbb{N}_0^d $ such that
	 the Taylor coeffcients, $c_{\boldsymbol{\alpha}}(f_n)) = Q_{\boldsymbol{\alpha}}\left( (c_{\boldsymbol{\beta}}(g_n))_{\boldsymbol{\beta} \leq \boldsymbol{\alpha}} \right)$.
	But $c_{\boldsymbol{\beta}}(g_n)) = c_{\boldsymbol{\beta}}(g))$ for $\boldsymbol{\beta} \leq \boldsymbol{n}$. 
	 Since, the rational functions $Q_{\boldsymbol{\alpha}}$ are universal, $Q_{\boldsymbol{\alpha}}\left( (c_{\boldsymbol{\beta}}(g_n))_{\boldsymbol{\beta} \leq \boldsymbol{\alpha}} \right)$ is equal to $c_{\boldsymbol{\alpha}}( \frac{g-1}{g+1})= c_{\boldsymbol{\alpha}}(f)$ for each $0 \leq \boldsymbol{\alpha} \leq \boldsymbol{n}$.
		
	Thus, the rational inner functions $f_n$ and $f$ both have the same Taylor polynomial up to order $\boldsymbol{n}$ for every $n$. Also, the multi-degree of the Cayley rational inner function $g_n$ is atmost $(2(n+1)^d, \dots, 2(n+1)^d)$ and hence it the multi-degree of $f_n$ is also atmost 
	$(2(n+1)^d, \dots, 2(n+1)^d)$ for every $n$. 
	
	Since, the family of rational inner functions $\{f_n : n \in \mathbb{N} \}$ is uniformly bounded by one on $\D^d$ and Taylor polynomials of $f_n$ and $f$ are the same up to order $\boldsymbol{n}$, $f_n$ converges to $f$ uniformly on compact subsets of $\D^d$ as $n$ tends to infinity. This completes the proof.
\end{proof}

\subsection{Pfister's approximation scheme}
In general, a bounded holomorphic function $f$ in the polydisk $\D^d$ can be approximated in the Fr\'echet topology of ${\mathcal O}(\D^d)$ by a sequence of polynomials $p_n$. Pfister proposes in his dissertation \cite{Pfister-thesis}
the following construction to approximate $f \in \mathcal{S}(\D^d) $. Pass first to a homothety $f_\rho(\boldsymbol{z}) = f(\rho \boldsymbol{z})$ with $0 < \rho <1$. The function $f_\rho$ can be uniformly approximated on the closed polydisk by a sequence of polynomials $p_n$.
Assume that the total degree of $p_n$ is $\kappa(n)$. Then its multi-degree is at most $N(n) = (\kappa(n), \kappa(n), \ldots, \kappa(n))$. Denote by $p_n^\ast$ the reflection of the polynomial $p_n$ with respect to the multi-index $N(n)$:
$$ p_n^\ast(\boldsymbol{z}) = \boldsymbol{z}^{N(n)}  \overline{p_n \left(\frac{1}{ \bar{\boldsymbol{z}}}\right)}, \ \ n \geq 1.$$
Then the sequence of rational inner functions
$$\phi_n(\boldsymbol{z}) = \frac{p_n(\boldsymbol{z}) + z_1 \ldots z_d \boldsymbol{z}^{N(n)}}{1 + z_1\ldots z_d p_n^\ast(\boldsymbol{z})}, \ \ n \geq 1,$$
converges to $f_\rho(\boldsymbol{z})$ in the topology of ${\mathcal O}(\D^d)$. Moreover, the Taylor polynomial of total degree $\kappa(n)$ of $\phi(\boldsymbol{z})$ coincides with $p_n(\boldsymbol{z})$. See Satz 14 in \cite{Pfister-thesis}. This simple approximation result was also independently discovered by Rudin and Stout \cite{Rudin-Stout}. In a recent work \cite{Alpay-BLMS}, this result for matrix-valued functions was established on the bidisk using an operator-theoretic tool, namely, unitary dilation of a contraction.  

\begin{remark}
Two notable differences between Pfister's approximation for functions in $\mathcal{S}(\D^d)$ and our Theorem \ref{CF-approx} for functions in $\mathcal{A}\mathcal{S}(\D^d)$ are: in the later the Taylor coefficients of the original function are matched increasingly, step by step, and second, the degrees of the rational approximation are not controlled in the former. 

Moreover, the proof of Theorem \ref{CF-approx} will extend to matrix-valued functions seamlessly as Theorem \ref{Thm:Inner-solution} can also be proved for matrix-valued interpolation data.
\end{remark}

\subsection{Pad\'e type approximation }

We propose a Pad\'e type approximation procedure via a sequence of anti-linear eigenvalue problems originating in the pioneering works of Takagi \cite{Takagi} and Pfister \cite{Pfister}. In its turn, this is a familiar Hankel form
approximation adapted to the multi-degree filtration in the space of trigonometric polynomials.

Let $f: \D^d \rightarrow \overline{\D}$ be a holomorphic function. For a fixed multi-index $\boldsymbol{n}$ in $\mathbb{N}_{0}^d$, we define two operators $T_{\boldsymbol{n}}$ and $C_{\boldsymbol{n}}$ on the finite dimensional space $\mathbb{C}_{\boldsymbol{n}}[\boldsymbol{z}]$ given by 
$$
T_{\boldsymbol{n}}(q) = P_{\boldsymbol{n}} (fq)  \ \text{ and } \ C_{\boldsymbol{n}}(q)(\boldsymbol{z}) = q^*(\boldsymbol{z}) =  \boldsymbol{z}^{\boldsymbol{n}} \overline{q\left(\frac{1}{\bar{\boldsymbol{z}}}\right)}, 
$$
for $q$ in $\mathbb{C}_{\boldsymbol{n}}[\boldsymbol{z}]$. Note that, $T_{\boldsymbol{n}}$ is linear but $C_{\boldsymbol{n}}$ is anti-linear.

\vspace{2mm}
\noindent \textbf{$C_{\boldsymbol{n}}$-symmetry of $T_{\boldsymbol{n}}$.} Note that $C_{\boldsymbol{n}}$ is involutive ($C_{\boldsymbol{n}}^2 = I_{\mathbb{C}_{\boldsymbol{n}}[\boldsymbol{z}]} $) and isomteric i.e., $\langle p, q \rangle = \langle C_{\boldsymbol{n}}(q), C_{\boldsymbol{n}}(p)\rangle$ for $p, q \in \mathbb{C}_{\boldsymbol{n}}[\boldsymbol{z}]$. Then $T_{\boldsymbol{n}}$ is $C_{\boldsymbol{n}}$-symmetric because for $p, q \in \mathbb{C}_{\boldsymbol{n}}[\boldsymbol{z}]$,
\begin{align*}
\left \langle C_{\boldsymbol{n}} T_{\boldsymbol{n}}  C_{\boldsymbol{n}}(p), q  \right \rangle =  \left \langle C_{\boldsymbol{n}} (q), T_{\boldsymbol{n}}  C_{\boldsymbol{n}}(p)  \right \rangle =  \left \langle q^*, fp^*  \right \rangle = \langle \bar{f} p, q \rangle = \left \langle T_{\boldsymbol{n}}^* p, q \right \rangle  
\end{align*}
implying $C_{\boldsymbol{n}} T_{\boldsymbol{n}}  C_{\boldsymbol{n}} =  T_{\boldsymbol{n}}^* $. See \cite{GP-TAMS07} for more details on $C$-symmetric operators. 

\subsection*{Takagi-Pfister anti-linear eigenvalue problem}
For the trancated Toeplitz operator $T_{\boldsymbol{n}}$ which is also $C_{\boldsymbol{n}}$-symmetric, the {\em Takagi-Pfister anti-linear eigenvalue problem} identifies non-negative eigenvalues $\sigma \in \mathbb{C}$ and non-zero polynomial eigenfunctions $q \in \mathbb{C}_{\boldsymbol{n}}[\boldsymbol{z}]$ such that 
\begin{align} \label{Eig-problem}
	T_{\boldsymbol{n} }C_{\boldsymbol{n}}(q) = \sigma q.
\end{align}
The anti-linearity of the first term of the equation allows to change the phase of the eigenvalue, so that always $\sigma \geq 0$.
Note the symmetry relation
$$ T_{\boldsymbol{n}} C_{\boldsymbol{n}}(q) = \sigma q, \ \ T^\ast_{\boldsymbol{n}} q = \sigma C_{\boldsymbol{n}}(q),$$
which gives rise to a Schmidt pair $(q, C_{\boldsymbol{n}}(q))$ for the operator $T_{\boldsymbol{n}}$.

The operator norm of $T_{\boldsymbol{n}}$ coincides with the highest eigenvalue:
$$
\| T_{\boldsymbol{n}} \| = \sigma_{\boldsymbol{n}} = \sup \{ \sigma \geq 0: T_{\boldsymbol{n}}  C_{\boldsymbol{n}} q = \sigma q , \ \exists q \neq 0 \}.
$$
For a proof see \cite[Corollary 4]{GP-TAMS07}.
Notice that $ \sigma_{\boldsymbol{n}} \leq \| f \|_\infty$. Denoting $\min {\boldsymbol{n}} = \min \{ n_1, n_2, \ldots, n_d\}$ one finds
$$ \lim_{\min  {\boldsymbol{n}} \rightarrow \infty}  \sigma_{\boldsymbol{n}} = \lim_{\min  {\boldsymbol{n}} \rightarrow \infty} \| T_ {\boldsymbol{n}} \| = \| f \|_\infty.$$

\subsection*{Decay of the remainder}
We carry on the notation and assumptions of the last section.
One can write $f C_{\boldsymbol{n}}q_{\boldsymbol{n}} = \sigma_{\boldsymbol{n}} q_{\boldsymbol{n}} + r_{\boldsymbol{n}}$
and seek an estimate of the decay of the remainder $r_{\boldsymbol{n}}$:
\begin{align} \label{Eq: decay}
	\| r_{\boldsymbol{n}}\|_2  \leq (\|f \|_\infty -\sigma_{\boldsymbol{n}}) \| q_{\boldsymbol{n}} \|_2.
\end{align}
Regarded as a Fourier or power series, $r_{\boldsymbol{n}}$ contains only terms involving monomials $z^{\boldsymbol{\alpha}}$ with $\boldsymbol{\alpha} \nleq \boldsymbol{n}$:
\begin{align*}
	r_{\boldsymbol{n}} (\boldsymbol{z}) = \sum_{\boldsymbol{\alpha} \nleq \boldsymbol{n}} c_{\boldsymbol{\alpha}} \boldsymbol{z}^{\boldsymbol{\alpha}}.
\end{align*}
Consider any $ 0 < \delta <1$ and  $|z_j| \leq \delta $ for $ 1\leq j \leq d$. Let $m= \min \boldsymbol{n}$.  Then
Cauchy-Schwarz inequality yields
\begin{align*}
	|r_{\boldsymbol{n}}(\boldsymbol{z})|  & \leq  \left( \sum_{\boldsymbol{\alpha} \nleq \boldsymbol{n}} |c_{\boldsymbol{\alpha}}|^2 \right)^{\frac{1}{2}} \left(\sum_{\boldsymbol{\alpha} \nleq \boldsymbol{n}} |\boldsymbol{z}^{\boldsymbol{\alpha}}|^2 \right)^{\frac{1}{2}} 
\end{align*}
where the sum in the second term in the right hand side is finite as $\boldsymbol{z} \in \delta \D^d$. Therefore,
\begin{align} 
	|r_{\boldsymbol{n}}(\boldsymbol{z})|  & \leq  (\|f \|_\infty - \sigma_{\boldsymbol{n}}) \| q_{\boldsymbol{n}} \|_2  \left(  \sum_{\boldsymbol{\alpha} \nleq \boldsymbol{n}}  \delta^{2|\boldsymbol{\alpha}|}  \right)^{\frac{1}{2}}  \notag \\
	& \leq  (\|f \|_\infty -\sigma_{\boldsymbol{n}}) \| q_{\boldsymbol{n}} \|_2  \left(  \sum_{|\boldsymbol{\beta}| \geq m}   \delta^{2|\boldsymbol{\beta}|} \right)^{\frac{1}{2}}  \notag \\
	& \leq (\|f \|_\infty -\sigma_{\boldsymbol{n}}) \| q_{\boldsymbol{n}} \|_2  \frac{\delta^m}{(1-\delta^2)^{d/2}}. \label{Decay estimate}
\end{align}

\subsection*{Detecting rational inner functions}
Suppose $f : \D^d \rightarrow \overline{\D}$ is holomorphic and for some $\boldsymbol{n} \in \mathbb{N}_0^d$, the largest positive eigenvalue $\sigma_{\boldsymbol{n}}$ of the problem \eqref{Eig-problem} is $1$ with normalized eigenfunction $q_{\boldsymbol{n}} \in \mathbb{C}_{\boldsymbol{n}}[\boldsymbol{z}]$. Then $\|f\|_\infty =1$ and the equation \eqref{Eq: decay} shows that $r_{\boldsymbol{n}} =0$. Thus
$f C_{\boldsymbol{n}}q_{\boldsymbol{n}} = q_{\boldsymbol{n}} $ and so, $f = \frac{q_{\boldsymbol{n}}}{C_{\boldsymbol{n}}q_{\boldsymbol{n}}}$ except possibly on the zero set of the polynomial $q_{\boldsymbol{n}}$ in $\D^d$. But, $f$ being bounded in $\D^d$, the rational function $ \frac{q_{\boldsymbol{n}}}{C_{\boldsymbol{n}}q_{\boldsymbol{n}}}$ must have holomorphic extension in $\D^d$ which means that $ q_{\boldsymbol{n}}$ is non-vanishing in $\D^d$. This proves that $f = \frac{q_{\boldsymbol{n}}}{C_{\boldsymbol{n}}q_{\boldsymbol{n}}}$ is a rational inner function on $\D^d$. We should also note that, in this case any eigenfunction $q_{\boldsymbol{n}} $ corresponding to $\sigma_{\boldsymbol{n}}=1$ of the eigenvalue problem \eqref{Eig-problem} produces the same rational function which is equal to $f$.

We warn the reader on a fundamental dichotomy between rational approximation by inner functions and polynomial approximation. Specifically,  the Takagi-Pfister approximation scheme does not finish in finitely many steps on a polynomial which is not a monomial. Indeed, Schwarz reflection implies that a polynomial depending on a single complex variable which has modulus equal to one on the unit circle is a monomial.

\subsection*{The Hankel form}
Let $f: \D^d \rightarrow \overline{\D}$ be a holomorphic function and fix a multi-index $\boldsymbol{n} \in \mathbb{N}_0^d$. The largest positive eigenvalue of the Takagi-Pfister problem \eqref{Eig-problem} for $2\boldsymbol{n}$ can be obtained from the optimization problem:
\begin{align} \label{Eq:optimum}
\text{Maximize}, \ \Re	\left \langle  T_{2 \boldsymbol{n} }C_{2 \boldsymbol{n}}(q) , q  \right \rangle = \Re \int_{\mathbb{T}^d} {f(\boldsymbol{\zeta})} \boldsymbol{\zeta}^{2\boldsymbol{n}}   \overline{q(\boldsymbol{\zeta})}^2 dm(\boldsymbol{\zeta}),
\end{align}
where $q \in \mathbb{C}_{2 \boldsymbol{n}}[\boldsymbol{z}]$ with $\| q\|_2 =1$. Note that the collection of all trigonometric polynomials of the form $\boldsymbol{\zeta}^{\boldsymbol{n}}   \overline{q(\boldsymbol{\zeta})}$ for $q \in \mathbb{C}_{2 \boldsymbol{n}}[\boldsymbol{z}]$, is the same as the collection $\operatorname{Trig}_{\boldsymbol{n}}$ of trigonometric polynomials whose Fourier coefficients vanish outside the rectangle $\prod_{j=1}^{d} [-n_j, n_j] \cap \mathbb{Z}^d$. We denote by $\operatorname{Trig}$ the vector space of all trigonometric polynomials.
The optimization problem \eqref{Eq:optimum} can be restated as 
\begin{align} \label{Eq:optimum-2}
	\text{Maximize}, \  \Re \int_{\mathbb{T}^d} {f(\boldsymbol{\zeta})}   q(\boldsymbol{\zeta})^2 dm(\boldsymbol{\zeta}),
\end{align}
where $q \in \operatorname{Trig}_{\boldsymbol{n}}$ with $\| q\|_2 =1$. 

This prompts to define the bounded Hankel form on $\operatorname{Trig}$ associated to the symbol $f$:
$$
H_f (q, r) = \langle f q , \bar{r} \rangle \ \text{ for } q,r \in \operatorname{Trig}.
$$
Thus the optimization problem \eqref{Eq:optimum-2} is equivalent to finding the maximum of $\Re H_f(q, q)$ with normalized $q$ varying over $\operatorname{Trig}_{\boldsymbol{n}}$. If the optimal value $\sigma_{2\boldsymbol{n}}$ is attained at $q_{2\boldsymbol{n}}$ in $\operatorname{Trig}_{\boldsymbol{n}}$, then the polynomial $\boldsymbol{z}^{\boldsymbol{n}} q_{2\boldsymbol{n}}$ in $\mathbb{C}_{\boldsymbol{n}}[\boldsymbol{z}]$ is an eigenfunction of the Takagi-Pfister problem \eqref{Eq:optimum}. We record the following simple observation whose proof follows from \eqref{Eq:optimum-2} and the fact that $\sigma_{2\boldsymbol{n}}$ converges to $\| f\|_{\infty}$ as $\min{\boldsymbol{n}}$ goes to infinity.

\begin{proposition} Let $f: \D^d \rightarrow \overline{\D}$ be a holomorphic function, continuous on the closed polydisk. Denote for every multi-index ${\boldsymbol{n}}$ an optimal solution $q_{2\boldsymbol{n}}$ of
(\ref{Eq:optimum-2}). Any weak-star limit point of the probability measures $|q_{2 \boldsymbol{n}}(\boldsymbol{\zeta})|^2 dm(\boldsymbol{\zeta})$ as $\lim \min{\boldsymbol{n}} = \infty$ is a positive measure supported on the set
$$ \{ \boldsymbol{\zeta} \in \mathbb{T}^d : \ |f(\boldsymbol{\zeta})| = \| f \|_\infty \}.$$
\end{proposition}

The Pad\'e type approximation scheme is the following. Given a holomorphic function $f: \D^d \rightarrow \overline{\D}$, consider a sequence of multi-indices $\boldsymbol{n}$ such that $\min \boldsymbol{n}$ tends to infinity.
Let $\sigma_{\boldsymbol{n}}$ be the highest eigenvalue of Takagi-Pfister problem, with associated eigenvector $q_{\boldsymbol{n}}$ such that $\| q_{\boldsymbol{n}} \| =1.$ Denote for symplicity $q^\ast_{\boldsymbol{n}} = C_{\boldsymbol{n}}
q_{\boldsymbol{n}}.$ The identity
$$ f q^\ast_{\boldsymbol{n}} = \sigma_{\boldsymbol{n}} q_{\boldsymbol{n}} + r_{\boldsymbol{n}}$$
can be interpreted as a formal series matching 
$$ f =  \sigma_{\boldsymbol{n}} \frac{q_{\boldsymbol{n}}}{ q^\ast_{\boldsymbol{n}}} + \frac{ r_{\boldsymbol{n}}}{ q^\ast_{\boldsymbol{n}}}$$
between $f$ and the rational function  $\sigma_{\boldsymbol{n}} \frac{q_{\boldsymbol{n}}}{ q^\ast_{\boldsymbol{n}}}$.
Indeed, if $q^\ast_{\boldsymbol{n}}(0) \neq 0$, then the germs of power series at $\boldsymbol{z}=0$ of $f$ and $\sigma_{\boldsymbol{n}} \frac{q_{\boldsymbol{n}}}{ q^\ast_{\boldsymbol{n}}}$ coincide up to the multi-index $\boldsymbol{n}$.
As expected, the location of the poles of these rational approximants is essential for an asymptotic evaluation of the convergence.

\begin{example}
	Consider the Schur function $f(z,w) = \frac{z+w}{2}$ on $\D^2$. If we take $\boldsymbol{n}=(1,1) \in \mathbb{N}^2$, then the eigenfunction of the Takagi-Pfister anti-linear eigenvalue problem corresponding to the largest eigenvalue $\sigma_{(1,1)} = \frac{1}{\sqrt{2}}$ is $p(z,w)= \frac{1}{2}(z+w)+ \frac{1}{\sqrt{2}} zw $. Therefore, $p^*(z, w) = \frac{1}{2}(z+w)+ \frac{1}{\sqrt{2}} $. Note that, the zero set of the polynomial $p^*$ intersects $\D^2$ non-trivially. Moreover, $p$ and $p^*$ are realtively prime.
\end{example}

\subsection*{The zero set problem}
Unlike the Schur-Takagi problem \cite{Takagi} in the disk, the polynomial $q^*$ can have zeros in $\D^d$ for an eigenfunction $q$ of the Takagi-Pfister problem \eqref{Eig-problem} corresponding to the largest positive eigenvalue, as evident from the above example.

If for each $\boldsymbol{n}\in \mathbb{N}^d$, we can ensure that there is an eigenfunction $q_{\boldsymbol{n}}$ corresponding to the largest positive eigenvalue of the problem \eqref{Eig-problem} such that the rational functions $\frac{q_{\boldsymbol{n}}}{q_{\boldsymbol{n}}^*}$ do not have poles in $\D^d$, then by a normal family argument one can prove that these rational functions, in fact rational inner functions $\frac{q_{\boldsymbol{n}}}{q_{\boldsymbol{n}}^*}$ converges uniformly on compact subsets of $\D^d$. 

The following example shows that, if the Schur function is a tensor product of one variable functions, then $\frac{q_{\boldsymbol{n}}}{q_{\boldsymbol{n}}^*}$ defines rational inner functions with poles off $\overline{\D^2}$.
\begin{example}
	Suppose $g$ and $h$ are two Schur functions on $\D$. Consider the Schur function $f(z, w)= g(z) h(w)$ on $\D^2$. For $\boldsymbol{n}=(n_1, n_2)$, we look at the one variable version of the problem \eqref{Eig-problem} for the functions $g$ and $h$.  Let $\lambda$ and $\sigma$ be the largest positive eigenvalue of the anti-linear eigenvalue problem of degree $n_1$ and $n_2$ related to $g$ and $h$ respectively. Let $p$ and $q$ be two eigenfunctions such that
	$$
	g(z) z^{n_1} \overline{p(1/\bar{z})} - \lambda p(z) = R(z) \ \text{ and }  h(w) w^{n_2} \overline{q(1/\bar{w})} - \sigma q(w) = S(w)
	$$
with the remainders, $R(z)=o(z^{n_1 +1})$ and $S(w)=o(w^{n_2 +1})$.	Then 
	\begin{align*}
		&f(z, w) z^{n_1} \overline{p(1/\bar{z})} w^{n_2} \overline{q(1/\bar{w})} - \lambda \sigma p(z)q(w) \\
		=& \ h(w) w^{n_2} \overline{q(1/\bar{w})} \left( g(z) z^{n_1} \overline{p(1/\bar{z})} - \lambda p(z)  \right) + \\
	&	\lambda p(z) S(w) \left(   h(w) w^{n_2} \overline{q(1/\bar{w})} - \sigma q(w)  \right) \\
		=& \ h(w) w^{n_2} R(z) + \lambda p(z) S(w)
	\end{align*}
	which implies 
	$$
	T_{\boldsymbol{n}} C_{\boldsymbol{n}} (p(z) q(w)) = \lambda \sigma p(z) q(w).
	$$
	Moreover, $	T_{\boldsymbol{n}}  C_{\boldsymbol{n}} = T_{n_1}^g  C_{n_1} \otimes T_{n_2}^h  C_{n_2} $ where 
	$$
	T_{n_1}^g C_{n_1} = \begin{bmatrix}
		0 & \dots & 0 & g_0 \\  0 & \dots & g_0 & g_1 \\ \vdots & \iddots& \vdots & \vdots \\ g_0 & \dots & g_{n_1-1} & g_{n_1}
	\end{bmatrix} \text{ and } T_{n_2}^h C_{n_2} = \begin{bmatrix}
	0 & \dots & 0 & h_0 \\  0 & \dots & h_0 & h_1 \\ \vdots & \iddots& \vdots & \vdots \\ h_0 & \dots & h_{n_2-1} & h_{n_2}
	\end{bmatrix}.
	$$
	Since, $\lambda$ and $\sigma$ are the largest eigenvalues of $T_{n_1}^g C_{n_1} $ and $T_{n_2}^h C_{n_2} $, respectively, the largest eigenvalue of $T_{\boldsymbol{n}} C_{\boldsymbol{n}}$ is $\lambda \sigma$.
	Finally, from the one variable result \cite[Theorem I]{Takagi}, we conclude that the poles of the rational function $\frac{p(z) q(w)}{C_{\boldsymbol{n}} (p(z) q(w))}$ are outside $\overline{\D^2}$.
\end{example}

In general, ensuring that the poles $\frac{q_{\boldsymbol{n}}}{q_{\boldsymbol{n}}^*}$ are disjoint from $\D^d$ seems to be a difficult task. Similar unresolved complications arise in multivariate digital filter/optimal polynomial approximation questions \cite{Beneteau, GW, Sargent-Sola} and in a variety of multi-variate Pad\'e approximation problems \cite{Cuyt, Guillaume-Huard}. It is natural at this moment to relax the hypothesis by adding some asymptotic assumptions, as for instance in the following theorem.

\begin{theorem}\label{weak-conv}
Let $f: \D^d \rightarrow \overline{\D}$ be a holomorphic function normalized to $\| f \|_\infty =1$ and let, for every multi-index ${\boldsymbol{n}},$ $q_{2\boldsymbol{n}}$ be an optimal solution of $L^2$-norm one of
(\ref{Eq:optimum-2}). Assume that the identically zero function is not a weak limit point of the sequence of complex analytic polynomials, $z^{\boldsymbol{n}} \overline{q_{2\boldsymbol{n}}}$ in the Hardy space $H^2(\D^d)$, as $\lim \min{\boldsymbol{n}} = \infty$.

Then there exists an analytic complex hypersurface $X \subset \D^d$ and a subsequence $\boldsymbol{n}(k), \  k \geq 1,\ \lim_k  \min \boldsymbol{n}(k) = \infty,$ such that the poles of the rational functions having unimodular values almost everywhere on $\mathbb{T}^d$,
$$ R_k(\boldsymbol{z}) = \frac{ q_{2 \boldsymbol{n}(k)} }{\overline{q_{2 \boldsymbol{n}(k)}}}, \ \ k \geq 1,$$
are accumulating to $X$ and $R_k$ converges to $f$ uniformly on every compact subset of $\D^d\setminus X$.
\end{theorem}

A clarification of the terminology: we say that the poles of the rational functions $R_k$ are accumulating to $X$, if, for every compact subset $K \subset \D^d \setminus X$ there exists a positive integer $M_K$ with the property that the functions $R_k$ are holomorphic in a neighborhood of $K$ for every $k \geq M_K$. Note that the ambient space is the open polydisk.

\begin{proof}

First we check that for every multi-index $\boldsymbol{n}$ and a trigonometric polynomial $p \in \operatorname{Trig}_{\boldsymbol{n}}$, the fraction $p/\overline{p}$ defines a rational function.
Indeed, $h(\boldsymbol{z}) = z^{\boldsymbol{n}} p(\boldsymbol{z})$ and $h^\ast(\boldsymbol{z}) =  \boldsymbol{z}^{\boldsymbol{n}} \overline{p(\boldsymbol{z})}$ are complex analytic polynomials and $h^\ast = C_{2\boldsymbol{n}} h$. Consequently, for the solution
$q_{2\boldsymbol{n}}$ of the Hankel optimization problem (\ref{Eq:optimum-2}) one identifies $h_{2 \boldsymbol{n}} = \boldsymbol{z}^{\boldsymbol{n}} q_{2\boldsymbol{n}}$ and  $h^\ast_{2 \boldsymbol{n}} = 
\boldsymbol{z}^{\boldsymbol{n}} \overline{q_{2\boldsymbol{n}}}$ as the optimal elements in Takagi-Pfister eigenvalue problem
$$ f h^\ast_{2 \boldsymbol{n}} = \sigma_{2 \boldsymbol{n}} h_{2 \boldsymbol{n}} + r_{2 \boldsymbol{n}}.$$

In view of the hypothesis, there exists a subsequence of extremal eigenfunctions $h^\ast_{2 \boldsymbol{n}(k)}$ converging weakly to a non-zero function $g \in H^2(\D^d)$. By passing to another subsequence we infer from
a normal family argument that $h^\ast_{2 \boldsymbol{n}(k)}$ converges to $g$ in the topology of $\mathcal{O}(\D^d)$. Let $X$ denote the zero set of $g$. 

Fix a compact subset $K$ of $\D^d \setminus X$. The rational functions $R_k = \frac {h_{2 \boldsymbol{n}(k)}}{h^\ast_{2 \boldsymbol{n}(k)}}$  are eventually analytic in a neighborhood of $K$.
Then, for these values of the index $k$:
$$ f(\boldsymbol{z}) = \sigma_{2\boldsymbol{n}(k)} R_k(\boldsymbol{z}) + \frac{ r_{2\boldsymbol{n}(k)}}{h^\ast_{2 \boldsymbol{n}(k)}}.$$
Since $g$ is non-vanishing in $K$, there is a $k_0 \in \mathbb{N}$ such that $\inf_{K} |h^\ast_{2 \boldsymbol{n}(k)}| >0$ for $k \geq k_0$.
Again $\sigma_{2\boldsymbol{n}(k)} \rightarrow \| f \|_\infty = 1$ and the remainder $r_{2\boldsymbol{n}(k)}$ converges uniformly to zero on $K$ as proved in (\ref{Decay estimate}).

In conclusion, $\lim_k \| f - R_k \|_{\infty, K} = 0.$
\end{proof}

\section{Final remarks and open questions} \label{Final remarks and open questions}

We collect a few open questions and related comments arising from this article.

\subsection*{The coefficient problem}
We proved that the set $K_{\boldsymbol{n}}$ of truncated Taylor series coefficients of a function belonging to Agler's class on the polydisk is semi-algebraic. Is $K_{\boldsymbol{n}}$ a basic semi-algebraic set,
 that is defined by a system of polynomial equalities and inequalities (without taking unions)? General criteria in this respect are known \cite{ABV}. Is it true in general for all Herglotz-Nevanlinna functions in the polydisk, or a
 classical domain in $\C^d, d >2,$ that the set of coefficients of the truncated power series expansion at a point is semi-algebraic? 

\subsection*{Convexity matters} The boundary of $K_{\boldsymbol{n}}$ in Euclidean space is also semi-algebraic. Are the coefficients of Cayley inner functions dense in it? Pfister proves by a simple example that $K_{\boldsymbol{n}}$
is not equal to the convex hull of the subset corresponding to Cayley inner functions. See Section 6 in \cite{Pfister-thesis}. A great deal of structural properties of convex semi-algebraic sets were recently discovered \cite{Sinn}.
Is there a linear matrix inequality description of $K_{\boldsymbol{n}}$? The single variable case suggests so \cite{Ahiezer-Krein}. We refer to \cite{NP} for
 a recent account of the fascinating topics of linear matrix inequality representations of convex semi-algebraic sets.

\subsection*{Convergence of the Pad\'e approximation} The zero distribution of Takagi-Pfister extremal functions $h^\ast_{\boldsymbol{n}}$ is critical for the convergence of the Pad\'e approximation scheme. This asymptotic analysis question remains wide open.

\end{document}